\documentclass[10pt,draftcls,onecolumn]{IEEEtran}

\usepackage{color}
\usepackage{enumerate}
\usepackage{epsfig}
\usepackage{multicol}
\usepackage{multirow}
\usepackage{amsmath,amsfonts,amsbsy,amssymb}
\usepackage{graphicx}
\usepackage{arydshln}  % For the dashed vertical line in matrix/array
\usepackage[subnum]{cases} % Sub-numbering equations within array
\usepackage{subfig}  % subfigures
\usepackage{cite}  % subfigures
\usepackage{algorithmic} % IEEE recommended algorithm style, or another package
\usepackage{threeparttable} % IEEE recommended footnote-in-table style (footnote)
\usepackage{mathtools} % for the correct symbol of definition, := \coloneqq
\usepackage{pifont} % for the dont-care term in the table
\usepackage{listings} % for code fragments into texts
\usepackage{lscape} % to landsscape/sideway/rotate long table

\allowdisplaybreaks[1]  % aligned multi-page equations
\newcommand{\beq}{\begin{eqnarray}}
\newcommand{\beqa}{\begin{array}}
\newcommand{\eeqa}{\end{array}}
\newcommand{\eeq}{\end{eqnarray}}
\newcommand{\beqstar}{\begin{eqnarray*}}
\newcommand{\eeqstar}{\end{eqnarray*}}

\newcommand{\bfa}{{\bf a}}
\newcommand{\bfb}{{\bf b}}
\newcommand{\bfc}{{\bf c}}
\newcommand{\bfx}{{\bf x}}

\newcommand{\bfq}{{\bf q}}

\newcommand{\bfy}{{\bf y}}

\newcommand{\bfw}{{\bf w}}

\newcommand{\bfe}{{\bf e}}

\newcommand{\bfk}{{\bf k}}
\newcommand{\bfz}{{\bf z}}

\newcommand{\bfv}{{\bf v}}

\newcommand{\bfepsilon}{{\boldsymbol \epsilon}}

\newcommand{\bfiota}{{\boldsymbol \iota}}

\newcommand{\bfrho}{{\boldsymbol \rho}}

\newcommand{\bftau}{{\boldsymbol \tau}}

\newcommand{\bfvarepsilon}{{\boldsymbol \varepsilon}}
\newcommand{\bfvarphi}{{\boldsymbol \varphi}}

\newcommand{\bfzeta}{{\boldsymbol \zeta}}

\newcommand{\bfzero}{{\bf 0}}
\newcommand{\real}{\mathbb{R}}

\newcommand{\calF}{{\mathcal F}}

\newcommand{\calI}{{\mathcal I}}

\newcommand{\calL}{{\mathcal L}}

\newcommand{\calN}{{\mathcal N}}

\newcommand{\calR}{{\mathcal R}}

\newcommand{\calX}{{\mathcal X}}

\newcommand{\iti}{{\it i}}

\newcommand{\itk}{{\it k}}
\newcommand{\itm}{{\it m}}
\newcommand{\itn}{{\it n}}
\newcommand{\itr}{{\it r}}

\newcommand{\comment}[1]{}

\newtheorem{Algorithm}{\indent\it Algorithm} % IEEE algorithm style as above
\newtheorem{Corollary}{\indent\it Corollary}

\newtheorem{Example}{\indent\it Example}
\newtheorem{Lemma}{\indent\it Lemma}
\newtheorem{Problem}{\indent\it Problem}
\newtheorem{Remark}{\indent\it Remark}
\newtheorem{Theorem}{\indent\it Theorem}

\newenvironment{enumerate-i-a}{
  \renewcommand{\theenumi}{\roman{enumi}}     % level-1 \ref:  i, ii, ...
  \renewcommand{\theenumii}{)\alph{enumii}}   % level-2 \ref:  i)a...
  \begin{enumerate}
}{
  \end{enumerate}
}

\newenvironment{enumerate-a}{
  \renewcommand{\theenumi}{\alph{enumi}}      % level-1 \ref:  a, b, c...
  \begin{enumerate}
}{
  \end{enumerate}
}

\newenvironment{enumerate-A-a}{
  \renewcommand{\theenumi}{\Alph{enumi}}      % level-1 \ref:  A, B, C...
  \renewcommand{\theenumii}{)\alph{enumii}}   % level-2 \ref:  A)a...
  \begin{enumerate}
}{
  \end{enumerate}
}
\begin{document}

\title{Exact Optimization: Part I}

\author{Li-Gang~Lin*~and~Yew-Wen~Liang% <-this % stops a space

\thanks{L.-G.~Lin (*corresponding author) is with the Department of Mechanical Engineering, National Central University, Taoyuan, 32001 Taiwan, and was with the Department of Computational Methods in Systems and Control Theory, Max Planck Institute for Dynamics of Complex Technical Systems, 39106 Magdeburg, Germany. (e-mail: lglin@g.ncu.edu.tw).}% <-this % stops a space
\thanks{Y.-W.~Liang is with the Institute of Electrical Control Engineering, National Yang Ming Chiao Tung University, Hsinchu 30010, Taiwan (e-mail: ywliang@cn.nctu.edu.tw).}% <-this stops a space
}

\maketitle

% According to TJM website, max: 150 words
\begin{abstract}\vspace{-0.3cm}
Nonlinear programming is explicitly analyzed via a novel perspective/method and from a bottom-up manner. The philosophy is based on the recent findings on convex quadratic equation (CQE), which help clarify a geometric interpretation that relates CQE to convex quadratic function (CQF). More specifically, regarding the solvability of CQE, its necessary and sufficient condition as well as a unified parameterization of all the solutions has recently been analytically formulated. Moving forward, the understanding of CQE is utilized to describe the geometric structure of CQF, and the CQE-CQF relation. All these results are shown closely related to a basis in the optimization literature, namely quadratic programming (QP). For the first time from this viewpoint, the QPs subject to equality, inequality, equality-and-inequality, and extended constraints can be algebraically solved in derivative-free closed formulae, respectively. All the results are derived without knowing a feasible point, a priori and any time during the process.
\end{abstract}

\vspace{-0.1cm}\def\abstractname{Key words and phrases}
\begin{abstract}
\noindent\vspace{-0.9cm}
\begin{center}
Nonlinear/quadratic programming, parametric optimization, matrix algebra, convex quadratic function/equation
\end{center}
\end{abstract}

\def\abstractname{2020 Mathematics Subject Classification}
\begin{abstract}
\noindent\vspace{-0.9cm}\begin{center}
  15A18, 52A41, 90C20, 90C25, 90C46
\end{center}
\end{abstract}

\section{Introduction}
\label{Sec_Intro}
Recently, \cite{LiLiHs:20} provided more understanding of the mathematical fundamentals: convex quadratic equation and function (CQE and CQF). In particular, \cite{LiLiHs:20} has analytically and completely formulated an equivalent solvability condition of CQE and a parameterization of all its solutions. These results encourage new investigations into a spectrum of applications, such as the one to nonlinear optimal control (\cite[Sec.~4]{LiLiHs:20}). The findings not only facilitate preliminary optimality recovery using the state-dependent Riccati equation (SDRE) scheme and its differential variant (see \cite{Cim:10} for a survey), but also enhance the computational performance of the other major application in this article. Specifically, we analyze the nonlinear programming \cite{Lu(Ye):03(16)}/convex optimization \cite{BoVa:04} from a new perspective based on \cite{LiLiHs:20}, which particularly boasts impacts in management science, operations research, and control engineering, for example, the model predictive control (MPC) \cite{CiFa:12}. It is worth mentioning that \cite{KaWaLiShPe:18} and \cite{QuSu:18} reveal more connection between the two applications, that is, nonlinear MPC using ``state-dependent coefficient'' (a design flexibility) in the SDRE scheme. Such a unified framework is envisioned/pioneered by \cite{BeMoDuPi:02}, and also pursued by this article.

With regard to the convex optimization, a basis in the literature of nonlinear programming algorithms is the quadratic programming (QP) \cite{Lu(Ye):03(16)}, which includes the linear programming (LP) as a special case and is closely linked with the CQE and CQF \cite{PaSa:10,RaBoJoVeJo:08}. Regarding the former reference, \cite{PaSa:10} utilize a solution subset of CQE, which is reformulated from the explicitly considered CQF and in terms of the polyhedral characterization \cite{Man:88}. Considering the latter, \cite{RaBoJoVeJo:08} make use of an important equality that is related with a CQF, as well as the existence and uniqueness of the corresponding optimum (point), with respect to a QP-like optimization problem. This supports the construction of a feasible sequence converging to the optimum. Generally speaking, more iterative/numerical QP solvers are available in the mainstream literature \cite{NoWr:06}, which are mostly analyzed from the differential perspective and can be classified by three levels \cite{BoVa:04}: i) equality-constrained QP; ii) linear equality-constrained optimization problem with twice-differentiable objective, as solved using Newton's method by reducing it to a sequence of equality-constrained QPs; and, after further imposing inequality constraints, iii) interior-point methods, which reduce the problem to a sequence of ii). It is worth noting that the equality-constrained QP solver is the most fundamental to build up various algorithms, subject to the constraints at the three levels. In addition, one area of great creativity among existing algorithms, namely the conjugate direction methods, illustrates that \textit{detailed analysis of QP leads to significant practical advances} \cite{Lu(Ye):03(16)}, such as solution accuracy, computation time, and more reported in \cite{JoFoTo:05} (see Sec. \ref{Subsec_Literature_Comparison} later for a comparison discussion). Moreover, the constrained QP can also be regarded as an optimization problem with the objective of constrained convex quadratic function (CCQF). These observations quite endorse the importance of QP and CQF, while motivate this article to provide further theoretical support (extended from \cite{LiLiHs:20}). \textit{One highlight is the first closed-form QP solver subject to a variety of constraints}. The analytical philosophy is inspired by \cite{Lue:69} (vector space methods) as well as the computational interest among a variety of research fields \cite{VaXiYa:20,DuJoWaWi:15,JoFo:13,Lu(Ye):03(16),RaLe:19}; but based on a novel understanding and perspective of the associated CQE, and interpretations of its relations to CQF and QP, respectively.

Following the previous findings of CQE and CQF \cite{LiLiHs:20}, this article deals with another application to nonlinear programming and -- in response to, for example, the open question in \cite[Sec. VI]{DuJoWaWi:15} -- the main result is exactly, autonomously, and algebraically solving the QPs that
\begin{itemize}
\item are subject to equality, inequality, equality-and-inequality, and extended constraints, respectively;
\item refrain from any derivative (only functional information is used);
\item consist of finite steps, the maximum of which is explicitly pre-determined;
\item and do not require any knowledge of a feasible point (any time before and during the process).
\end{itemize}
The derivation leverages the vector space methods such as the singular value decomposition (SVD), and investigates from a perspective based on \cite{LiLiHs:20} and geometric interpretations that clarify the structures of CQE, CQF, and CQE-CQF relation, respectively.

The rest of this article is arranged as follows. The notations and the formulation of main problem are provided in Sec. \ref{Sec_Notation_Problem}. The method proposed in \cite{LiLiHs:20} is applied to this focused application ``nonlinear programming/convex optimization'' in Sec. \ref{Sec_Convex_Optimization}. At first, the analysis of extreme/minimal values of CQF is performed in Sec. \ref{Subsec_Unconstrained_QP}, which is followed by a complete characterization of the explicit solution to the equality-constrained QP in Sec. \ref{Subsec_Equality_QP}. As a step forward, the two preliminary subsections are shown indispensable for the proposed closed-form QP solver in Sec. \ref{Subsec_QP}, and an extension subject to the constraints in Sec. \ref{Subsec_Extended_QP}. Therefore, discussions on differences and advantages with respect to existing solvers in the literature can be justified in Sec. \ref{Subsec_Literature_Comparison}. Finally, three representative examples are demonstrated to offer computational evidences and additional insights (as verified using MATLAB$\textsuperscript{\textregistered}$) in Sec. \ref{Sec_Ex}, while concluding remarks are given in Sec. \ref{Sec_Conclusion}.

\section{Notation and Problem Formulation}
\label{Sec_Notation_Problem}
We use the notational conventions, unless noted otherwise: the symbols $||\cdot||$, $\calN(\cdot),~\calR(\cdot),~\lambda(\cdot),~(\cdot)^\dagger,~(\cdot)^{\dagger/2}$, and $(\cdot)^T$ denote the Euclidean norm, null space, range space, spectrum/eigenvalues, Moore-Penrose generalized inverse (shortened as ``pseudoinverse'' afterward), square root of the pseudoinverse \cite{UrBlDr:13}, and transpose of a vector or matrix, respectively. Additionally, we denote $(\cdot)^\perp$ as the orthogonal complement of a vector space, $\oplus$ the direct sum of vector spaces, $\bfe_1$ the first standard basis vector in $\real^\itn$, and $\{\sigma_i\}_{i=1}^\itn$ the set of singular values of an $\itn\times\itn$ matrix with rank $\itr$ ($\sigma_1\ge\sigma_2\ge\cdots\ge\sigma_\itn\ge 0$ and $\sigma_\itr>0$). Moreover, denote $P\succ 0$ (resp., $P\succeq 0$), if $P=P^T\in\real^{{\it n}\times {\it n}}$ is positive definite (resp., semidefinite) \cite{Lu(Ye):03(16)}.

Consider the following CQF $F_\bfx: \real^\itn\rightarrow\real$ \cite{Lu(Ye):03(16)},
\beq
F_\bfx(\bfx)=\bfx^TP\bfx/2+\bfq^T\bfx+s,
\label{CQF}
\eeq
where both $\bfx,~\bfq\in\real^{\it n}$, $P=P^T\in\real^{{\it n}\times {\it n}}$, $P\succeq 0$, and $s\in\real$. {In particular, the function (\ref{CQF}) is strictly convex, if its Hessian matrix $P\succ 0$}. Note that the early findings in \cite{LiLiHs:20} start with the formulation of CQE, as in Eq. (\ref{CQE}) later, to emphasize its dominance. Then, in this continuing work, the concept of CQF dominates, and thus we also give its formulation in the very beginning (\ref{CQF}) and interpret its close relation to CQE later in detail. More specifically, \textit{in this article, the main focus shifts to explicitly analyzing and solving the following problem.}

\begin{Problem}
From a new viewpoint based on the solution to {\rm\cite[Problem~2.1]{LiLiHs:20}} as well as the CQE-CQF relation, solve the QPs in an analytical, autonomous, and derivative-free manner and subject to equality, inequality, equality-and-inequality, and extended constraints, respectively.
\label{Prob_Convex_Optimization}
\end{Problem}

\section{Application to Convex Optimization (QP)}
\label{Sec_Convex_Optimization}
\textit{The main concepts} behind all the derivations in this section are the geometrical interpretations of i) the value/difference ``$\bfq^TP^\dagger\bfq/2-s$'' according to the solvability conditions in Eqs. (\ref{Solvability_Cond_n}) and (\ref{Solvability_Cond_k_R(M)}) later and ii) the solution parameterization in Eqs. (\ref{Solution_n}) and (\ref{Solution_k_R(M)}). Notably, in i), $P^\dagger=P^{-1}$ when $P$ is nonsingular. Moreover, the novelty can be revealed by that i) facilitates the derivations that exploit the (additional) \textit{analysis perspective} from the image of the CQF (\ref{CQF}). As a step forward, the \textit{hierarchical layers} in the parameterization, with respect to ii), exactly categorize the location(s) of the unconstrained optimum/optima. In the following, Sec. \ref{Subsec_Unconstrained_QP} considers the unconstrained optimization, which are preliminaries supporting the derivations of constrained ones in the remaining subsections. Following the QP formulation/definition in literature \cite{BoVa:04,Lu(Ye):03(16)}, Sec. \ref{Subsec_Extended_QP} enlarges the included types of constraints, without introducing any excessive variable but, actually, reducing to a lower-dimensional unconstrained equivalent problem. Finally, Sec. \ref{Subsec_Literature_Comparison} discusses advantages or differences over existing solvers \cite{NoWr:06}. Remarkably, we recall the following lemma because it is largely used in this article, where the notation is consistent with \cite{LiLiHs:20,LiLiHs:draft} to avoid misleading/confusion.

\begin{Lemma} (Solvability and Solutions of CQE {\rm\cite[Theorem 3.1]{LiLiHs:20}})\\
Consider the CQE below {\rm\cite{Lu(Ye):03(16)}}:
\beq
\bfz^TM\bfz+\bfk^T\bfz+c=0,
\label{CQE}
\eeq
where both $\bfz,~\bfk\in\mathbb{R}^{\it n}$, $M=M^T\in\mathbb{R}^{{\it n}\times {\it n}}$, $M\succeq 0$, and $c\in\mathbb{R}$.
\begin{enumerate-A-a}
\item\label{Lem_Solutions_n} If rank$(M)=n$, then CQE {\rm(\ref{CQE})} is solvable, if and only if (iff)
    \beq
    \bfk^TM^{-1}\bfk\ge 4c.
    \label{Solvability_Cond_n}
    \eeq
    Accordingly, the set of solutions are, and can be parameterized by,
    \beq
    \bfz=-M^{-1}\bfk/2+\sqrt{\bfk^TM^{-1}\bfk/4-c}\cdot M^{-1/2}\cdot \bfv,
    \label{Solution_n}
    \eeq
    where $\bfv\in\mathbb{R}^{\it n}$ and $||\bfv||=1$.
\item\label{Lem_Solutions_r} Otherwise (rank$(M)<n$), it is solvable, iff {\rm(\ref{Solvability_Cond_k_R(M)})} or {\rm(\ref{Solvability_Cond_k_not_R(M)})}, where
    \beq
    &&\bfk\in \calR(M)~\mbox{and}~\bfk^TM^\dagger\bfk\ge 4c,\label{Solvability_Cond_k_R(M)}\\
    &&\bfk\not\in \calR(M).\label{Solvability_Cond_k_not_R(M)}
    \eeq
    Accordingly, the sets of solutions are, and can be parameterized by, respectively,
    \begin{enumerate}
    \item\label{Lem_Solutions_r_k_R(M)} for Condition {\rm(\ref{Solvability_Cond_k_R(M)})},
        \beq
        \bfz=-M^\dagger\bfk/2+\sqrt{\bfk^TM^\dagger \bfk/4-c}\cdot M^{\dagger/2}\bfrho+\bfvarepsilon,~~
        \label{Solution_k_R(M)}
        \eeq
        where both $\bfrho,~\bfvarepsilon\in\mathbb{R}^{\it n}$, $\bfrho\in \calR(M)$, $||\bfrho||=1$, and $\bfvarepsilon\in \calN(M)$;
    \item\label{Lem_Solutions_r_k_not_R(M)} for Condition {\rm(\ref{Solvability_Cond_k_not_R(M)})}, decompose $\bfk=\bfk_M+\bfk_{M^\perp}$, where $\bfk_M\in \calR(M)$, $\bfk_{M^\perp}\in \calR(M)^\perp$, and both $\bfk_M,~\bfk_{M^\perp}\in\mathbb{R}^\itn$. Then,
        \beq
        \bfz=-(F_\bfw/||\bfk_{M^\perp}||^2)\cdot\bfk_{M^\perp}+\bfvarphi+\bftau,
        \label{Solution_k_not_R(M)}
        \eeq
        where the CQF $F_\bfw: \calR(M)\subset\mathbb{R}^\itn\rightarrow\mathbb{R}$,
        \beq
        F_\bfw(\bfw)=\bfw^TM\bfw+\bfk_M^T\bfw+c,
        \label{CQF_Solution_k_not_R(M)}
        \eeq
        all $\bfw,~\bfvarphi,~\bftau\in\mathbb{R}^{\it n}$, both $\bfw,~\bftau\in \calR(M)$, and $\bfvarphi\in \calN(M)\cap \calN(\bfk^T)$.
    \end{enumerate}
\end{enumerate-A-a}
\label{Lem_Solutions}
\end{Lemma}

\subsection{Unconstrained QP}
\label{Subsec_Unconstrained_QP}

\begin{Theorem} (Solutions to Unconstrained QP)\\
Consider the optimization problem, minimize $F_\bfx$ {\rm(\ref{CQF})} with respect to $\bfx$.
\begin{enumerate-i-a}
\item\label{Thm_Unconstrained_QP_all_preimage} The preimage of any level set of $F_\bfx$ {\rm(\ref{CQF})} can be parameterized by Eqs. {{\rm(\ref{Solution_n})}, {\rm(\ref{Solution_k_R(M)})-(\ref{CQF_Solution_k_not_R(M)})}}, respectively, where $M=P/2$, $\bfk=\bfq$, $c=s-l$, and $l\in\real$ is any level set value of $F_\bfx$.
\item\label{Thm_Unconstrained_QP_finite_condition} The optimal value is finite, iff $~\bfq\in\calR(P)$.
\item\label{Thm_Unconstrained_QP_finite_value} The finite optimal value is $l^*=s-\bfq^T P^\dagger\bfq/2$. And the corresponding unique optimum is, or optima are all parameterized by, $\bfx^*=\bfx_p^*+\breve\bfvarepsilon$, where $\bfx_p^*\coloneqq-P^\dagger\bfq$, $\breve\bfvarepsilon\in\real^\itn$, and $\breve\bfvarepsilon\in \calN(P)$.
\end{enumerate-i-a}
\label{Thm_Unconstrained_QP}
\end{Theorem}

\begin{figure}[htbp]
    \begin{center}
    \includegraphics[width=5.8cm]{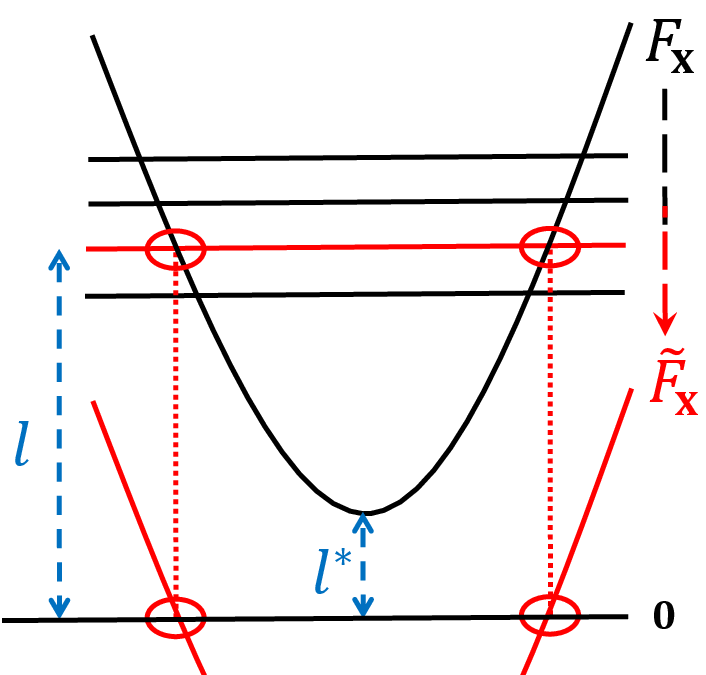}
    \caption{Geometric interpretation of a ``shift'' $l$ and the critical one ``$l^*$'' in Sec. {\rm\ref{Subsec_Unconstrained_QP}}.}
    \label{Fig_Shift}
    \end{center}
\end{figure}

\begin{figure}[htbp]
    \begin{center}
    \includegraphics[width=6.2cm]{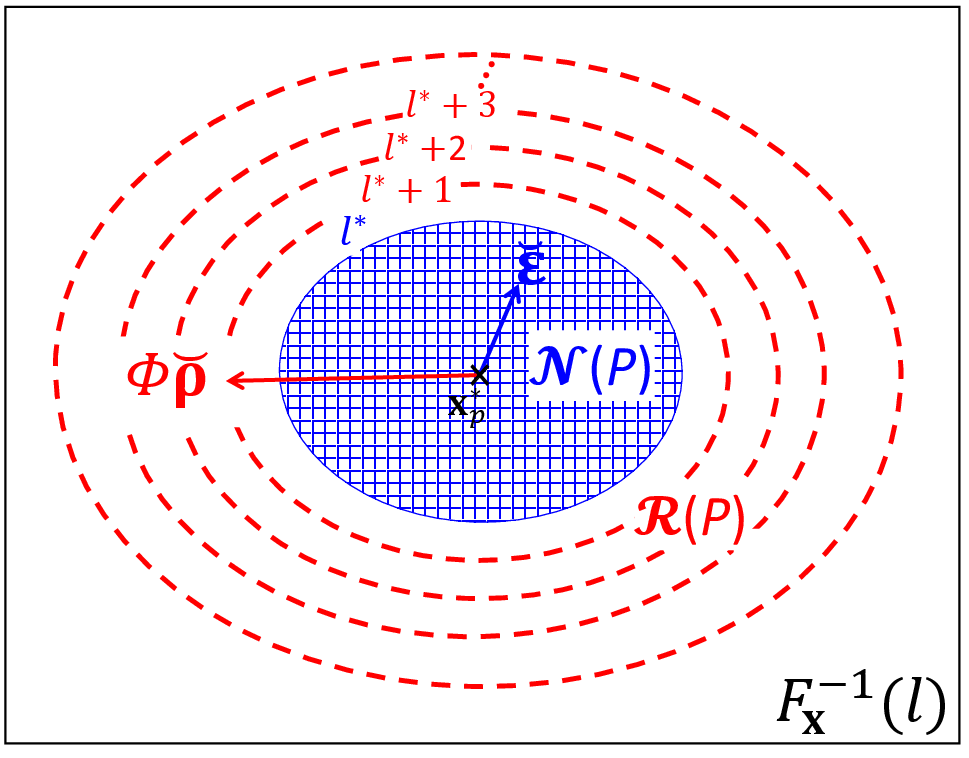}
    \caption{Geometric interpretation of ``hierarchical layers'' in Sec. {\rm\ref{Subsec_Unconstrained_QP}} (Remark {\rm\ref{Rem_Fig_Geometric_Layers}}).}
    \label{Fig_Layers}
    \end{center}
\end{figure}

\begin{IEEEproof}
See Appendix \ref{App_Proof_Unconstrained_QP} for a complete and unified proof, which is new with regard to \ref{Thm_Unconstrained_QP_finite_condition}) and \ref{Thm_Unconstrained_QP_finite_value}), while first for \ref{Thm_Unconstrained_QP_all_preimage}). Along with the proof, a geometric interpretation for the novel concept of ``critical shift'' (resp., ``hierarchical layers'') is given in Fig. \ref{Fig_Shift} (resp., Remark \ref{Rem_Fig_Geometric_Layers} and Fig. \ref{Fig_Layers}).
\end{IEEEproof}

\begin{Remark}
Based on Lemma \ref{Lem_Solutions}, the results of \ref{Thm_Unconstrained_QP_finite_condition}) and \ref{Thm_Unconstrained_QP_finite_value}) in Theorem \ref{Thm_Unconstrained_QP} are consistent with the literature, for example, \cite{BoVa:04} as analyzed from the differential perspective, which thus and also endorses Lemma \ref{Lem_Solutions}. All the results in Theorem \ref{Thm_Unconstrained_QP}, particularly \ref{Thm_Unconstrained_QP_all_preimage}), will be utilized for the constrained and extended QP problems, later in Secs. \ref{Subsec_Equality_QP} to \ref{Subsec_Extended_QP}. On the other hand, more potentials are expected toward the applications to least-squares approximation and regression analysis. Finally, according to the respective items in Theorem \ref{Thm_Unconstrained_QP}, it is worth remarking that:
\begin{enumerate-i-a}
\item A geometric interpretation is illustrated in Fig. \ref{Fig_Shift}.
\item From its proof, typically Eq. (\ref{Proof_Thm_Unconstrained_QP_CQF_n}), $F_\bfx$ is always unbounded above; bounded below, iff $\bfq\in\calR(P)$.
\item Regarding the special but popular case of nonsingular $P$ \cite{Lu(Ye):03(16)}, then $P^\dagger=P^{-1}$, $\calN(P)=\{\bfzero\}$, the unique optimum $\bfx^*=-P^{-1}\bfq$, and $l^*=s-\bfq^TP^{-1}\bfq/2$.
\end{enumerate-i-a}
\label{Rem_Unconstrained_QP}
\end{Remark}

\begin{Remark}
After introducing the concept of ``critical shift $l^*$'' in this section, as illustrated in Fig. \ref{Fig_Shift}, it suffices to complement the examination previously in \cite[Remark~3.2]{LiLiHs:20}. By \ref{Thm_Unconstrained_QP_finite_condition}) and \ref{Thm_Unconstrained_QP_finite_value}) of Theorem \ref{Thm_Unconstrained_QP}, the optimal/minimal value ($l^*$) is finite and equals zero (no shift) for $M\succeq 0$, either singular or nonsingular. In other words, $\bfz^TM\bfz>0$, for all $\bfz\not\in\calN(M)$. To sum up, combining \cite[Remark~3.2]{LiLiHs:20} and this Remark \ref{Rem_Thm_Unconstrained_QP_PSD} endorse both Theorem \ref{Thm_Unconstrained_QP} and Lemma \ref{Lem_Solutions} above, according to the basic mathematical properties of positive definiteness and semidefiniteness.
\label{Rem_Thm_Unconstrained_QP_PSD}
\end{Remark}

\subsection{Equality-Constrained QP}
\label{Subsec_Equality_QP}

\begin{Theorem} (Solutions to Equality-Constrained QP)\\
Consider the optimization problem,
\beq
\mbox{minimize}~F_\bfx~\mbox{subject to}~A\bfx=\bfb,
\label{Problem_Equality_QP}
\eeq
where $F_\bfx$ is in Eq. {\rm(\ref{CQF})}, $A\in\real^{\itm\times\itn}$, rank$(A)=\itm\in\{1,2,\cdots,\itn-1\}$, and $\bfb\in\real^\itm$. From the SVD of $A$, let a $V_2\in\real^{\itn\times(\itn-\itm)}$ be given, which is associated with $\calR(V_2)=\calN(A)$. Denote $\bar l\in\real$ as any level set value of $F_\bfx$.
\begin{enumerate-A-a}
\item\label{Thm_Equality_Constrained_QP_Condition_Unconstrained_QP} If $V_2^TPV_2\ne 0\Leftrightarrow\calR(P)\cap\calN(A)\ne\{\bfzero\}\Leftrightarrow\calN(A)\not\subseteq\calN(P)$,
    \begin{enumerate}
    \item\label{Thm_Equality_Constrained_QP_Unconstrained_QP_Qquivalence} Problem {\rm(\ref{Problem_Equality_QP})} is equivalent to an unconstrained QP,
    \item\label{Thm_Equality_Constrained_QP_Unconstrained_QP_All_Preimage} the preimage of any level set/value of ``$F_\bfx$ subject to $A\bfx=\bfb$'' can be respectively parameterized by
        \beq
        \bfx=A^\dagger\bfb+V_2\bfy,
        \label{Parameterization_Equality_Constraint}
        \eeq
        where $\bfy\in\real^{\itn-\itm}$ is determined by Eqs. {{\rm(\ref{Solution_n})}, {\rm(\ref{Solution_k_R(M)})}-(\ref{CQF_Solution_k_not_R(M)})}, with $\bfz=\bfy$, $M=V_2^TPV_2/2$, $\bfk=V_2^T(\bfq+PA^\dagger\bfb)$, and $~c=s+[\bfq^T+\bfb^T(A^\dagger)^TP/2]A^\dagger\bfb-\bar l$,
    \item\label{Thm_Equality_Constrained_QP_Unconstrained_QP_Finite_Condition} among all $\bar l$ in {\rm\ref{Thm_Equality_Constrained_QP_Unconstrained_QP_All_Preimage})}, the optimal one/minimum is finite, iff $~V_2^T(\bfq+PA^\dagger\bfb)\in\calR(V_2^TPV_2)$,
    \item\label{Thm_Equality_Constrained_QP_Unconstrained_QP_Optimal_Value_Optimum} following {\rm\ref{Thm_Equality_Constrained_QP_Unconstrained_QP_Finite_Condition})}, this finite optimal value is $\bar l^*\in\real$,
        \beq
        \bar l^*=[\bfb^T(A^\dagger)^TP/2+\bfq^T][I_\itn-V_2(V_2^TPV_2)^\dagger V_2^TP]A^\dagger\bfb+s-\bfq^TV_2(V_2^TPV_2)^\dagger V_2^T\bfq/2,
        \label{Optimal_Value_Equality_Constrained_QP_Unconstrained_QP}
        \eeq
        and the corresponding unique optimum is, or optima are all, parameterized by,
        \beq
        \bfx^{\bar *}=\bfx_p^{\bar *}+V_2\bfvarepsilon^{\bar *},
        \label{Optimum_Equality_Constrained_QP_Unconstrained_QP}
        \eeq
        where $\bfx_p^{\bar *}\coloneqq A^\dagger\bfb-V_2(V_2^TPV_2)^\dagger V_2^T(\bfq+PA^\dagger\bfb)\in\real^\itn$ denotes the particular solution of $\bfx^{\bar *}\in\real^\itn$, $\bfvarepsilon^{\bar *}\in\real^{\itn-\itm}$, and $\bfvarepsilon^{\bar *}\in\calN(V_2^TPV_2)$,
    \item\label{Thm_Equality_Constrained_QP_Unconstrained_QP_Uniqueness} to further categorize {\rm\ref{Thm_Equality_Constrained_QP_Unconstrained_QP_Optimal_Value_Optimum})}, the optimum is unique, iff $\calN(A)\cap\calN(P)=\{\bfzero\}$;
    \end{enumerate}
\item\label{Thm_Equality_Constrained_QP_Condition_Unconstrained_LP} else, if $V_2^T(\bfq+PA^\dagger\bfb)\ne\bfzero$,
    \begin{enumerate}
    \item Problem {\rm(\ref{Problem_Equality_QP})} is equivalent to an unconstrained LP, and thus unbounded,
    \item the preimage of any level set/value of ``$F_\bfx$ subject to $A\bfx=\bfb$'' can be respectively parameterized by Eq. (\ref{Parameterization_Equality_Constraint}), where $\bfy$ is determined by {\rm\cite[Lemma~3.1]{LiLiHs:20}}, with $\bfz=\bfy$, $\bfzeta=V_2^T(\bfq+PA^\dagger\bfb)$, and $\nu=s+[\bfq^T+\bfb^T(A^\dagger)^TP/2]A^\dagger\bfb-\bar l$;
    \end{enumerate}
\item\label{Thm_Equality_Constrained_QP_Condition_Constant} else,
    \begin{enumerate}
    \item Problem {\rm(\ref{Problem_Equality_QP})} is equivalent to a constant function,
    \item its value is ``$s+[\bfq^T+\bfb^T(A^\dagger)^TP/2]A^\dagger\bfb$'',
    \item the (only) preimage is parameterized by Eq. {\rm(\ref{Parameterization_Equality_Constraint})}, for all $\bfy$.
    \end{enumerate}
\end{enumerate-A-a}
\label{Thm_Equality_Constrained_QP}
\end{Theorem}

\begin{IEEEproof}
See Appendix \ref{App_Proof_Equality_Constrained_QP} for the first, complete, and unified proof.
\end{IEEEproof}

\begin{Remark}
Both Theorems \ref{Thm_Unconstrained_QP} and \ref{Thm_Equality_Constrained_QP} present the results within a unified framework. In particular, the following summarize further analyses of the more significant, quadratic case \ref{Thm_Equality_Constrained_QP_Condition_Unconstrained_QP}) in Theorem \ref{Thm_Equality_Constrained_QP}.
\begin{enumerate-i-a}
\item Using the common terminology, such as \cite[Sec. 4.1.1]{BoVa:04}, the considered problem is \text{feasible} since $\bfb\in\calR(A)$. The preimage of any level set of ``$F_\bfx$ (\ref{CQF}) subject to $A\bfx=\bfb$'' is the \textit{feasible set}. The definitions of the optimum/optima and corresponding optimal value follow the general consensus \cite{BoVa:04} throughout Sec. \ref{Sec_Convex_Optimization}, and thus omitted for brevity.
\item\label{Rem_Thm_Equality_Constrained_QP_Vertex} Regarding the consistent constraint $A\bfx=\bfb$ (in other words, $\bfb\in\calR(A)$ where $A$ is of full rank) in this theorem, we omit the only remaining case (that is, rank$(A)=\itn$) for a more concise presentation. This corresponds to $\calN(A)=\{\bfzero\}$ and that there is only one feasible point $\bfx^{\bar *}=A^\dagger\bfb$. However, by an extended definition of $V_2\coloneqq 0_\itn$, the results of {\ref{Thm_Equality_Constrained_QP_Unconstrained_QP_Finite_Condition})-\ref{Thm_Equality_Constrained_QP_Unconstrained_QP_Uniqueness})} in this theorem still apply to this sub-QP problem.
\item\label{Rem_Thm_Equality_Constrained_QP_Boundedness} In the proof, the considered QP Problem (\ref{Problem_Equality_QP}) is equivalently reformulated as an unconstrained QP: minimize $F_\bfy$ (\ref{Proof_Equality_Constrained_QP_Equivalent_y}) with respect to $\bfy$. This is viable, iff $V_2^TPV_2\ne 0$. With regard to the unconstrained QP, $F_\bfy$ (resp., the original QP Problem (\ref{Problem_Equality_QP})) is always unbounded above; bounded below, iff $V_2^T(\bfq+PA^\dagger\bfb)\in\calR(V_2^TPV_2)$. Note that a geometric interpretation regarding the CQF $F_\bfy$ (resp., QP (\ref{Problem_Equality_QP}) for the optimal value) can be inferred from Fig. \ref{Fig_Shift}.
\end{enumerate-i-a}
\label{Rem_Thm_Equality_Constrained_QP}
\end{Remark}

\begin{Remark}
An example to endorse Theorem \ref{Thm_Equality_Constrained_QP} comes from the minimum or least-norm problem \cite{BoVa:04,GoVa:13}. More specifically, to find the optimum $\bfx^{\bar *}_{ln}$ that (satisfies the consistent underdetermined constraint $A\bfx=\bfb$ and) is of least $l_2$-norm, which (norm) is most commonly adopted. This is in the applicable form using this theorem, with $P=2I_\itn\succ 0$, $\bfq=\bfzero$, and $s=0$. Accordingly, this is a QP given $\calR(P)\cap\calN(A)=\calN(A)$, which (null space) is at least one-dimensional and thus nonzero. Moreover, the optimal value (denoted $\bar l^*_{ln}$) is finite since $V_2^TA^\dagger\bfb=V_2^TV_1\Gamma_1^{-1}W^T\bfb=\bfzero\in\calR(I_{\itm})=\real^\itm$, where both $V_1\in\real^{\itn\times\itm}$ being orthogonal to $V_2$ and $A^\dagger=V_1\Gamma_1^{-1}W^T\in\real^{\itn\times\itm}$ are given in Eq. (\ref{App_Proof_Equality_Constrained_QP_A_dagger}) in the proof. Specifically,
\beq
\bar l^*_{ln}&=&\bfb^T(A^\dagger)^TA^\dagger\bfb-\bfb^T(A^\dagger)^TV_2(V_2^TV_2)^{-1} \underbrace{V_2^TA^\dagger}_\text{$0$}\bfb\nonumber\\
&=&\Vert A^\dagger\bfb\Vert^2.\nonumber
\eeq
Correspondingly and similarly, the unique optimum is
\beq
\bfx^{\bar *}_{ln}&=&A^\dagger\bfb-2V_2(V_2^TPV_2)^\dagger V_2^TA^\dagger\bfb+V_2\bfvarepsilon^{\bar *}\nonumber\\
&=&A^\dagger\bfb,
\eeq
where $\bfvarepsilon^{\bar *}\in\calN(V_2^TV_2)=(\real^{\itn-\itm})^\perp=\{\bfzero\}$. As expected, these results agree with the literature \cite{BoVa:04}.
\label{Rem_Least_Norm}
\end{Remark}

\begin{Remark}
This remark focuses upon the special case of $\itm=1$, and more efficiently constructs an example of $V_2$ in Theorem \ref{Thm_Equality_Constrained_QP}. Without loss of generality, let $A=\bfa^T\in\real^{1\times\itn}$ be normalized. Given the results in \cite[Theorem~4.2]{LiLiHs:20} and further analyses in \cite[Remark~4.4]{LiLiHs:20}, $V_2$ is readily available by computing the last $(\itn-1)$ columns of $H_\bfiota^\bfa\coloneqq I_\itn-2\bfiota_\bfa\bfiota_\bfa^T$, where $\bfiota_\bfa\coloneqq (\bfa-\bfe_1)/(\Vert\bfa-\bfe_1\Vert)$ and $Y=I_{\itn-1}$. This remark also indicates a research direction, which aims at extending to general $\itm$ in the equality-constrained QP (\ref{Problem_Equality_QP}) based on the generalization of \cite[Theorem 4.2]{LiLiHs:20}, and thus involves the orthogonalization on the rows of $A$ in the first place.
\label{Rem_Thm_Equality_Constrained_QP_m_1}
\end{Remark}

\subsection{QP}
\label{Subsec_QP}
Consider the optimization problem \cite{Lu(Ye):03(16)},
\beq
&&\mbox{minimize}~F_\bfx\nonumber\\
&&\mbox{subject to}~A\bfx=\bfb,\nonumber\\
&&\hspace{1.53cm}\bfc_i^T\bfx\le d_i,~i\in\calI=\{1,2,\cdots,\kappa\},
\label{Problem_QP}
\eeq
where $F_\bfx$ is in Eq. (\ref{CQF}), $A\in\real^{\itm\times\itn}$, rank$(A)=\itm\in\{1,2,\cdots,\itn-1\}$, $\bfb\in\real^\itm$, $\bfc_i\in\real^\itn$, $d_i\in\real$, and there exists a point that satisfies all the constraints.

\begin{figure}[htbp]
    \begin{center}
    \includegraphics[width=8.1cm]{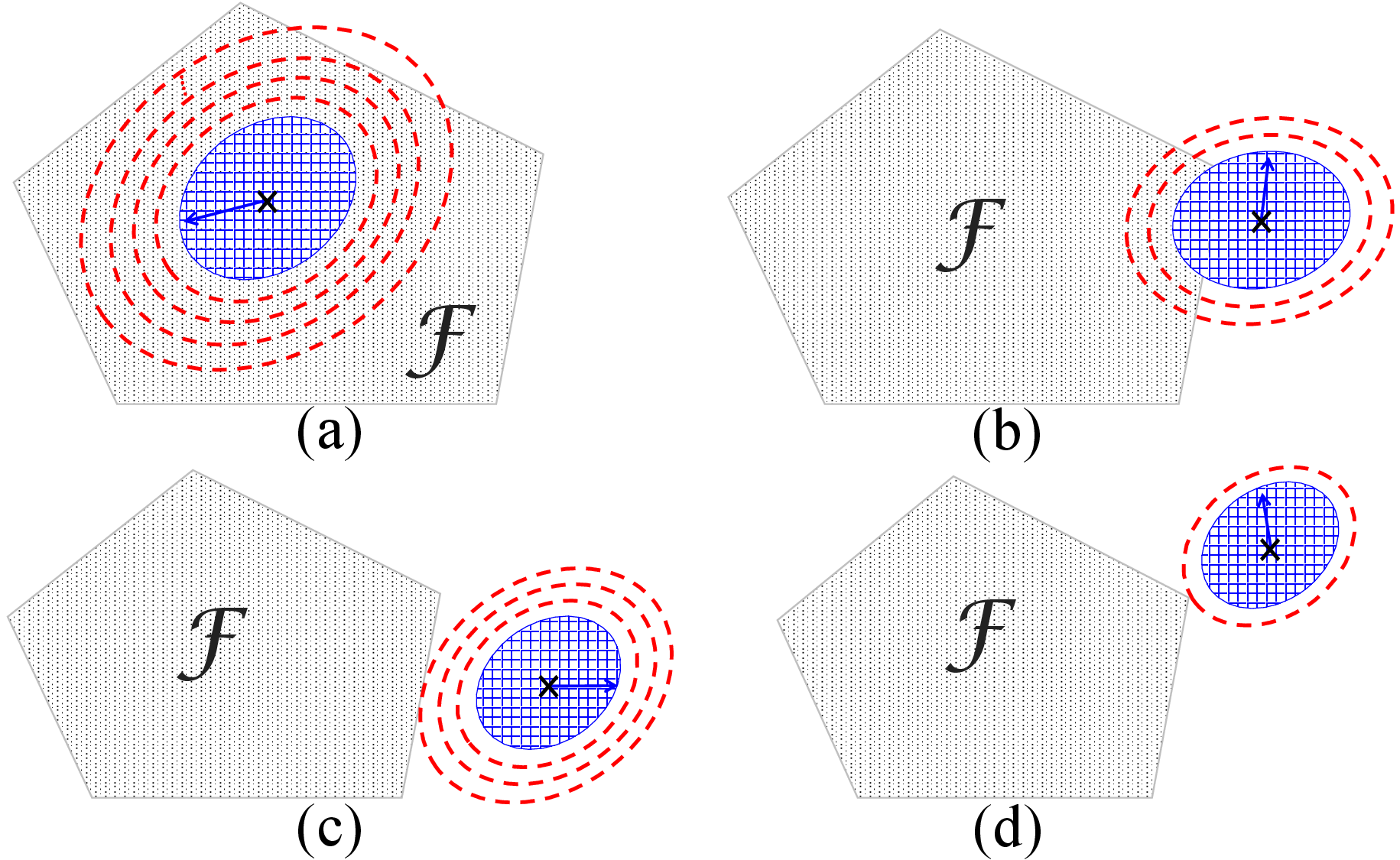}
    \caption{Representative cases in terms of the location of $\bfx^{\bar *}$ {\rm(\ref{Optimum_Equality_Constrained_QP_Unconstrained_QP})} for Algorithm {\rm\ref{Alg_QP}}: (a) $\bfx_p^{\bar *}\in\calF$, (b) $\bfx_p^{\bar *}\not\in\calF$ and $\bfx^{\bar *}\in\calF$ for some $\bfvarepsilon^{\bar *}\ne\bfzero$, (c) $\bfx^{\bar *}\not\in\calF$ while $\bfx^{\tilde *}$ on an edge, and (d) $\bfx^{\bar *}\not\in\calF$ while $\bfx^{\tilde *}$ on a vertex.}
    \label{Fig_Locations}
    \end{center}
\end{figure}

\begin{figure}[htbp]
    \begin{center}
    \includegraphics[width=8cm]{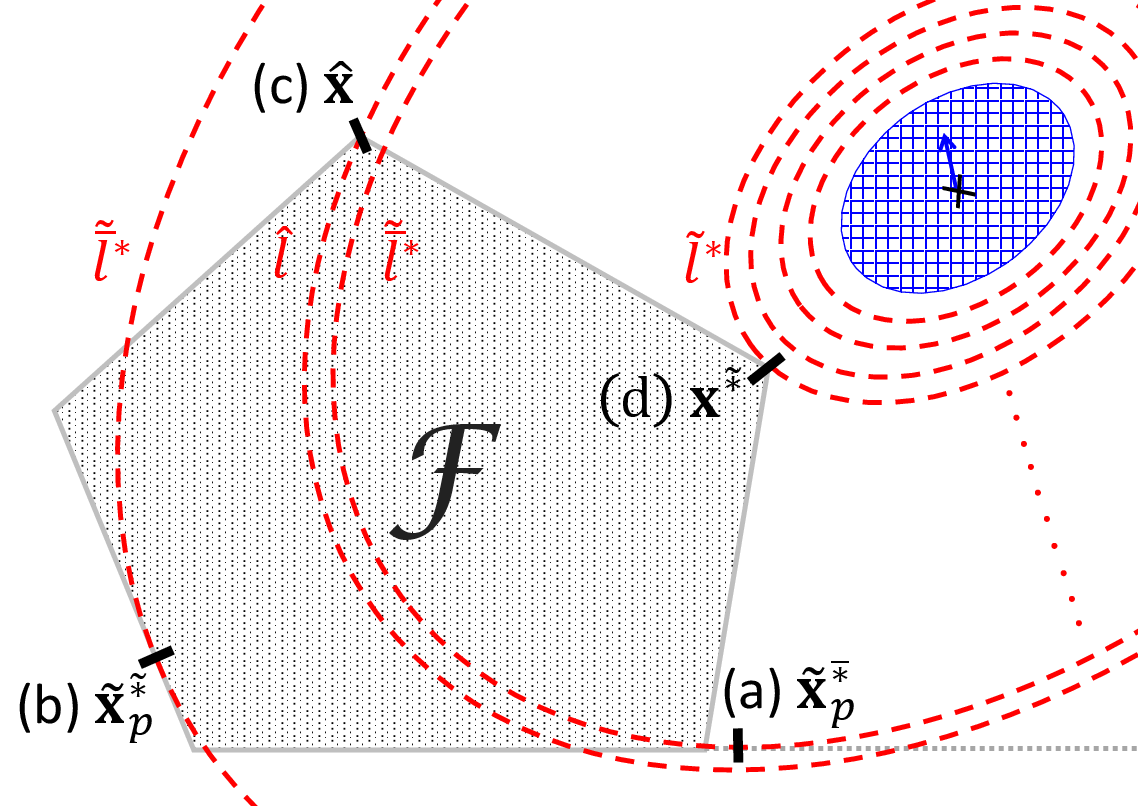}
    \caption{A planar demonstration of ``optimality candidates'' in Appendix {\rm\ref{App_Proof_QP_bcd}} when $\bfx^{\bar *}\not\in\calF$ according to Eq. {\rm(\ref{Optimum_Equality_Constrained_QP_Unconstrained_QP})}, where (a) $\tilde\bfx^{\bar *}$ is not a candidate, both (b) and (c) $\tilde\bfx^{\tilde *}$ are candidates on an edge and vertex, respectively, and the optimality occurs at another vertex/candidate (d) $\bfx^{\tilde *}$.}
    \label{Fig_Candidates}
    \end{center}
\end{figure}

To this general QP Problem (\ref{Problem_QP}), an explicit solver (namely, without any iterations, approximations, or assumptions as in the classical/numerical optimization \cite{NoWr:06,Lu(Ye):03(16)}) is available in Algorithm \ref{Alg_QP}, whose proof is in Appendix \ref{App_Proof_QP} (together with geometric demonstrations via Figs. \ref{Fig_Locations} and \ref{Fig_Candidates}).

\vspace{0.2cm}
\hrule width \hsize \kern 0.5mm \hrule width \hsize
\vspace{0.25cm}
\begin{Algorithm} Closed-Form QP Solver\\
\vspace{-0.45cm}
\hrule
\vspace{0.15cm}
\begin{algorithmic}[1]
\REQUIRE $P,~\bfq,~s,~A,~\bfb,~\{\bfc_i,d_i\}_{i\in\calI}$ as defined in {\rm(\ref{Problem_QP})}\\
\hspace{-0.59cm}\textbf{Initial Conditions:} $\bar \calL=\bar \calX=\bar \calI=\tilde \calI^*=\emptyset$
\ENSURE The optimal value $\tilde l^*$ and the (resp., a subset of) corresponding optimum (resp., optima) $\bfx^{\tilde *}$
\STATE\label{Alg_A_V_2} compute the SVD of $A$ till $\calR(V_2)=\calN(A)$ is obtained
\vspace{0.1cm}
\IF{$V_2^TPV_2\ne 0$}\vspace{0.1cm}\label{Alg_a_QP_Condition}
\STATE\label{Alg_a_optimum_particular} $\bfx_p^{\bar *}=A^\dagger\bfb-V_2(V_2^TPV_2)^\dagger V_2^T(\bfq+PA^\dagger\bfb)$
\STATE\label{Alg_a_optimal_value} $\bar l^*=[\bfb^T(A^\dagger)^TP/2+\bfq^T][I_\itn-V_2(V_2^TPV_2)^\dagger V_2^TP]A^\dagger\bfb+s-\bfq^TV_2(V_2^TPV_2)^\dagger V_2^T\bfq/2$\vspace{0.1cm}

\IF{``$V_2^TPV_2\succ 0$'' and ``$\bfc_i^T\bfx_p^{\bar *}\le d_i$ for all $i\in\calI$''}\vspace{0.1cm}\label{Alg_a_optimum_in}
\RETURN\label{Alg_Equality_QP_Positive_Definite_In} $\tilde l^*=\bar l^*$, $\bfx^{\tilde *}=\bfx_p^{\bar *}$, and $\tilde \calI^*$
\ELSIF{``$V_2^T(\bfq+PA^\dagger\bfb)\in\calR(V_2^TPV_2)$'' and ``$\bfc_i^T\bfx_p^{\bar *}\le d_i$ for all $i\in\calI~$''}\label{Alg_Equality_QP_Optimum_Condition_Inequality}
\vspace{0.1cm}
\STATE\label{Alg_Switch_Element} $(\bar l^*, \bfx_p^{\bar *},\emptyset)\in\bar\calL\times\bar\calX\times\bar\calI$
\STATE\label{Alg_Switch_Go_To} go to line {\rm\ref{Alg_For_All}}
\ENDIF\label{Alg_a_QP_Condition_End}

\ELSE\label{Alg_bcd}

\FORALL{$\calI_j\coloneqq\{j_1,j_2,\cdots,j_k\}\subseteq\calI$, where $k\in\{1,2,\cdots,\kappa\}$,}\label{Alg_For_All}

\STATE\label{Alg_bcd_Augmented_Parameters} let $\tilde C=[\bfc_{j_1},\bfc_{j_2},\cdots,\bfc_{j_k}]^T$, $\tilde A=[A^T,\tilde C^T]^T$, and $\tilde\bfb=[\bfb^T,d_{j_1},d_{j_2},\cdots,d_{j_k}]^T$\vspace{0.1cm}
\IF{$\tilde\bfb\in\calR(\tilde A)$}\label{Alg_bcd_Condition_Consistency}\vspace{0.1cm}
\IF{rank$(\tilde A)=\itm+\itk<\itn$}\vspace{0.1cm}\label{Alg_bcd_Condition_Full_Rank_Singular}
\STATE\label{Alg_bcd_A_V_2} compute the SVD of $\tilde A$ till $\calR(\tilde V_2)=\calN(\tilde A)$ is obtained\vspace{0.1cm}
\IF{``$\tilde V_2^TP\tilde V_2\ne 0$'',\vspace{0.1cm} ``$\tilde V_2^T(\bfq+P\tilde A^\dagger\tilde\bfb)\in\calR(\tilde V_2^TP\tilde V_2)$'', and ``$\bfc_i^T\tilde \bfx_p^{\bar *}\le d_i$ for all $i\in\calI\backslash \calI_j$,\vspace{0.1cm} where $\tilde\bfx_p^{\bar *}=\tilde A^\dagger\tilde\bfb-V_2(V_2^TPV_2)^\dagger V_2^T(\bfq+P\tilde A^\dagger\tilde\bfb)$''}\label{Alg_Inequality Equality_QP_Optimum_Condition_Inequality}
\STATE\label{Alg_bcd_Condition_Full_Rank_Singular_Optimal_Value} $\tilde{\bar{l}}^*=[\tilde\bfb^T(\tilde A^\dagger)^TP/2+\bfq^T][I_\itn-V_2(V_2^TPV_2)^\dagger V_2^TP]\tilde A^\dagger\tilde\bfb+s-\bfq^TV_2(V_2^TPV_2)^\dagger V_2^T\bfq/2$
\STATE\label{Alg_bcd_Condition_Full_Rank_Singular_Inclusion} $(\tilde{\bar{l}}^*,\tilde \bfx_p^{\bar *},\calI_j)\in\bar\calL\times\bar\calX\times\bar\calI$
\ENDIF\label{Alg_bcd_Condition_Full_Rank_Singular_End}

\ELSIF{rank$(\tilde A)=\itn$}\vspace{0.1cm}\label{Alg_bcd_Condition_Full_Column_Rank}

\STATE\label{Alg_bcd_Condition_Full_Column_Rank_Check_x} compute $\hat\bfx=\tilde A^\dagger\tilde\bfb$\vspace{0.1cm}
\IF{$\bfc_i^T\hat\bfx\le d_i$, for all $i\in\calI\backslash \calI_j$}\vspace{0.1cm}\label{Alg_bcd_Condition_Full_Column_Rank_Feasibility}
\STATE\label{Alg_bcd_Condition_Full_Column_Rank_Optimal_Value} $\hat l=F_\bfx(\hat\bfx)$, where $F_\bfx$ is from {\rm(\ref{CQF})}\vspace{0.1cm}
\STATE\label{Alg_bcd_Condition_Full_Column_Rank_Inclusion} $(\hat l,\hat\bfx,\calI_j)\in\bar\calL\times\bar\calX\times\bar\calI$\vspace{0.1cm}
\ENDIF\label{Alg_bcd_Condition_Full_Column_Rank_End}

\ENDIF\label{Alg_bcd_Condition_Full_Rank_End}
\ENDIF

\ENDFOR\vspace{0.1cm}\label{Alg_For_All_End}

\STATE\label{Alg_bcd_Result} $\tilde l^*=\min \bar\calL$, while $\bfx^{\tilde *}$ and $\tilde\calI^*$ are as associated\vspace{0.1cm}
\RETURN\label{Alg_bcd_Return} $\tilde l^*$, $\bfx^{\tilde *}$, and $\tilde \calI^*$
\ENDIF\label{Alg_End}
\end{algorithmic}
\label{Alg_QP}
\end{Algorithm}
\vspace{0.2cm}
\hrule width \hsize \kern 0.5mm \hrule width \hsize
\vspace{0.4cm}

\begin{Remark} Regarding Algorithm \ref{Alg_QP},
\begin{enumerate-i-a}
\item it does not require to know a feasible point \textit{at any time} in the solving process (in contrast to, for example, \textit{a priori} to initiate the process according to \cite{BoVa:04,Lu(Ye):03(16)}) since the core of the analytical philosophy is a novel categorization of QP Problem (\ref{Problem_QP}) in Appendix \ref{App_Proof_QP}, which particularly builds on new results in Lemma \ref{Lem_Solutions}.
\item\label{Rem_Alg_Order_For_All} in the for-all environment (lines \ref{Alg_For_All}-\ref{Alg_For_All_End}), the processing order is arbitrary, as long as every nonempty subset of $\calI$ is examined. Examples are given later in Sec. \ref{Sec_Ex}.
\item an economic/expedited version of this algorithm is available with respect to Case (a) in Fig. \ref{Fig_Locations}, as motivated by the particular interest in the computational performance using existing QP solvers (to name a few, \cite{JoFo:13,RaLe:19}). Specifically, by replacing lines \ref{Alg_Switch_Element} and \ref{Alg_Switch_Go_To} with \ref{Alg_Equality_QP_Positive_Definite_In}, the complete searching for all optima is economized since the optimal value and an optimum have already been obtained, without further effort for any possible, additional optimum. Accordingly, the algorithm is expedited because it terminates at line \ref{Alg_Switch_Element} without proceeding into the for-all environment (lines \ref{Alg_For_All}-\ref{Alg_For_All_End}). More details and related discussions can be found in the end of Appendix \ref{App_Proof_QP_a}, as well as Remark \ref{Rem_QP_Singleton_F}.
\end{enumerate-i-a}
\label{Rem_Alg}
\end{Remark}

\begin{Corollary} With regard to Algorithm {\rm\ref{Alg_QP}}, consider
\begin{enumerate-A-a}
\item\label{Cor_Alg_QP_P} a special case of the QP Problem {\rm(\ref{Problem_QP})}, $P\succ 0$ as is commonly seen in literature {\rm\cite{Lu(Ye):03(16)}}. Algorithm {\rm\ref{Alg_QP}} can be more efficiently compacted by removing lines {\rm\ref{Alg_a_QP_Condition}}, {\rm\ref{Alg_Equality_QP_Optimum_Condition_Inequality}}-{\rm\ref{Alg_Switch_Go_To}}, {\rm\ref{Alg_bcd}}, and {\rm\ref{Alg_End}}, as well as ``$V_2^TPV_2\succ 0$'' at line {\rm\ref{Alg_a_optimum_in}}, and both ``$\tilde V_2^TP\tilde V_2\ne 0$'' and ``$\tilde V_2^T(\bfq+P\tilde A^\dagger\tilde\bfb)\in\calR(\tilde V_2^TP\tilde V_2)$'' at line {\rm\ref{Alg_Inequality Equality_QP_Optimum_Condition_Inequality}}.
\item\label{Cor_Alg_Q_Inequality} an extension to the QP under inequality constraints only, that is, $A=0$ and $\bfb=\bfzero$ in the QP Problem {\rm(\ref{Problem_QP})}. Algorithm {\rm\ref{Alg_QP}} applies to this case by replacing ``both $A, A^\dagger$ with $0$'', ``$V_2$ with $I_\itn$'', ``$\bfx_p^{\bar *}$ with $\bfx_p^*$'', ``$\bfb$ with $\bfzero$'', and ``$\bar l^*$ with $l^*$''.
\end{enumerate-A-a}
\label{Cor_Alg_QP}
\end{Corollary}

\begin{IEEEproof}
\ref{Cor_Alg_QP_P}): Given $P\succ 0$, in Algorithm \ref{Alg_QP} we readily have $V_2^TPV_2\succ 0$, $V_2^T(\bfq+PA^\dagger\bfb)\in\calR(V_2^TPV_2)=\real^{\itn-\itm}$, $\tilde V_2^TP\tilde V_2\succ 0$, and $\tilde V_2^T(\bfq+P\tilde A^\dagger\tilde\bfb)\in\calR(\tilde V_2^TP\tilde V_2)=\real^{\itn-\itm-\itk}$. The result immediately follows.

\ref{Cor_Alg_Q_Inequality}): At first, the concept is to replace the application of Theorem \ref{Thm_Equality_Constrained_QP} by Theorem \ref{Thm_Unconstrained_QP} in the beginning of Algorithm \ref{Alg_QP}. More specifically, the processes before the for-all environment (lines \ref{Alg_a_optimum_particular}-\ref{Alg_a_QP_Condition_End}) instead consider the unconstrained QP, and determine whether the associated unique optimum, or the particular solution of nonunique optima, satisfies the (inequality) constraints. Starting from line \ref{Alg_For_All}, the remaining processes in Algorithm \ref{Alg_QP} are slightly modified. The main difference is that, at line \ref{Alg_bcd_Augmented_Parameters}, the remaining processes are instead built on the ``base'' of unconstrained QP ($A=0$ and $\bfb=\bfzero$) and then taking each inequality constraint into consideration. All in all, from this novel perspective and design, this extension can be easily established by the simple parameter replacements.
\end{IEEEproof}

\begin{Remark} With respect to Corollary \ref{Cor_Alg_QP},
\begin{enumerate-i-a}
\item\label{Rem_Cor_Alg_P_Nonsingular} in \ref{Cor_Alg_QP_P}) of Corollary \ref{Cor_Alg_QP}, the more compacted algorithm renders $\bfx^{\bar *}=\bfx_p^{\bar *}$ and $\tilde\bfx^{\bar *}=\tilde\bfx_p^{\bar *}$ for all $\calI_j\subseteq\calI$. That is, the associated solution freedom vanishes (because of $P\succ 0$). Even though Algorithm \ref{Alg_QP} is designed to \textit{implicitly} while \textit{exhaustively} include the optimality searching within such a freedom (more detailed arguments at Appendix \ref{App_Proof_QP}), this remark also reminds that the searching process, excessive in this case, automatically vanishes as well -- as designed.
\item Example \ref{Ex_QP_Nonsingular} in Sec. \ref{Sec_Ex} demonstrates both Algorithm \ref{Alg_QP} and Corollary \ref{Cor_Alg_QP}.
\end{enumerate-i-a}
\label{Rem_Cor_Alg}
\end{Remark}

\subsection{Extended QP}
\label{Subsec_Extended_QP}
\begin{Theorem} (Solutions to An Extended QP)\\
Consider the optimization problem,
\beq
\mbox{minimize}~F_\bfx~\mbox{subject to both}~\bfx,~\bfq\in\calR(P),
\label{Problem_Extended_QP}
\eeq
where $F_\bfx$ is in Eq. {\rm(\ref{CQF})}.
\begin{enumerate-A-a}
\item\label{Thm_Extended_QP_Preimage} The preimage of any level set of the CCQF ``$F_\bfx$ subject to both $\bfx,~\bfq\in\calR(P)$'' can be respectively parameterized by
    \beq
    \bfx=-P^{\dagger}\bfq+\sqrt{\bfq^TP^\dagger\bfq-2s+2\check l}\cdot P^{\dagger/2}\cdot \check \bfrho,~~~~~~~~~
    \label{Thm_Extended_QP_Solutions}
    \eeq
    where $\check l\in\real$ is any level set value of the CCQF, $\check\bfrho\in \calR(P)$, and $\Vert\check\bfrho\Vert=1$.
\item\label{Thm_Extended_QP_Optimality} The optimal value is finite, and equals $\check l^*=s-\bfq^TP^\dagger\bfq/2$. The corresponding unique optimum is $\bfx^{\hat *}=-P^\dagger\bfq$.
\end{enumerate-A-a}
\label{Thm_Extended_QP}
\end{Theorem}

\begin{IEEEproof}
See Appendix \ref{App_Proof_Extended_QP}.
\end{IEEEproof}

\begin{Corollary}
Consider the specific but popular case, $P\succ 0$ {\rm\cite{Lu(Ye):03(16)}}. Theorem {\rm\ref{Thm_Extended_QP}} can be specialized by replacing $P^\dagger$ with $P^{-1}$, $P^{\dagger/2}$ with $P^{-1/2}$, and ``$\check\bfrho\in\calR(P)$'' with ``$\check\bfrho\in\real^\itn$'', respectively.
\end{Corollary}

\begin{IEEEproof}
Given $P\succ 0$, such a case-specific result readily follows by relating the basic SVD properties in \cite[Eqs. {(35)-(39)}]{LiLiHs:20} to this special consideration.
\end{IEEEproof}

\begin{Remark}
An application of Theorem \ref{Thm_Extended_QP} is associated with the CQF
in Eq. (\ref{CQF_Solution_k_not_R(M)}), that is, to obtain its preimage $\bfw$ at any level set/value. To be more detailed, by letting $P=2M$ which is rank deficient, $\bfq=\bfk_M$, and $s=c$ in accordance with the CQF in Eq. (\ref{CQF}), the preimage $\bfw$ can be algebraically while completely parameterized by (using Theorem \ref{Thm_Extended_QP})
\beq
\bfw=-M^\dagger\bfk_M/2+\sqrt{\bfk_M^TM^\dagger\bfk_M/4-c+\breve F_\bfw}\cdot M^{\dagger/2}\cdot \hat\bfrho,
\label{Preimage_CQF_Solution_k_not_R(M)}
\eeq
where $\breve F_\bfw\in\real$ is any (given) level set value, $\hat\bfrho\in \calR(M)$, and $\Vert\hat\bfrho\Vert=1$. As an example, this application-benefit is also demonstrated (with an illustration) in \cite[Sec.~5]{LiLiHs:20}.
\label{Rem_Thm_Extended_QP_w}
\end{Remark}

\subsection{Comparisons and Impacts to Literature}
\label{Subsec_Literature_Comparison}
With respect to the QP Problem (\ref{Problem_QP}), several representative and popular solvers in the literature or readily in the market \cite{NoWr:06} are selected to demonstrate this section. At first, according to \cite{Lu(Ye):03(16)}, such a general QP (with inequality constraints) is usually, numerically solved using an ``Active Set Method'' (ASM), which is also valued in the fields of management science \cite{Hey:22,KoYa:91} and control engineering \cite{JoFo:13}, to name a few. It has been reported that a feasible point is required to initiate the solving process of ASM -- including a simpler/common case of $P\succ 0$ \cite{Lu(Ye):03(16)} -- till the latest MATLAB$\textsuperscript{\textregistered}$ implementation (namely, in the ``quadprog'' routine). Similar difficulty/inconvenience can also be found in other celebrated solvers: Newton's method, the barrier method, and primal-dual interior-point methods \cite{Hey:22,BoVa:04,JoFo:13}. \textit{In contrast}, the proposed solvers autonomously generate the results/solutions, notably Algorithm \ref{Alg_QP} and Theorem \ref{Thm_Equality_Constrained_QP}. The analytical philosophy differentiates from the classical Lagrange/Primal-Dual method \cite{BoVa:04,Lu(Ye):03(16),NoWr:06}, but stems from Lemma \ref{Lem_Solutions} as well as Figs. \ref{Fig_Shift} (``critical shift'') and \ref{Fig_Layers} (``hierarchical layers'').

In addition, another advantage is revealed from an accuracy perspective, with regard to the other representative solver: the barrier method that is renowned for large-scale problems \cite{JoFo:13}. The importance of such an interior-point solver is particularly focused upon in \cite[Chapter 11]{BoVa:04}. Accordingly, under reasonable assumptions, a minimum number of iteration steps are generally required for convergence, or a desired accuracy; whilst, the maximum number of steps (and computation time) are deemed challenging \cite{JoFo:13}. \textit{As a comparison}, the proposed results yield exact solutions for both the optimal point(s) and value. Specifically, considering the equality-constrained QP Problem (\ref{Problem_Equality_QP}), the finite solutions (when existed, according to a comprehensive categorization in Theorem \ref{Thm_Equality_Constrained_QP}) are available in Eqs. (\ref{Optimal_Value_Equality_Constrained_QP_Unconstrained_QP}) and (\ref{Optimum_Equality_Constrained_QP_Unconstrained_QP}). Moreover, regarding the general (resp., extended) QP problem in (\ref{Problem_QP}) (resp., (\ref{Problem_Extended_QP})), Algorithm \ref{Alg_QP} (resp., Theorem \ref{Thm_Extended_QP}) also renders results in analytical, algebraic, and closed forms. Remarkably, it is worth mentioning that all the proposed results refrain from any convergence analysis, namely the estimation of convergence rate, which can cause computational difficulty \cite{DuJoWaWi:15} but is generally necessary till the up-to-date literature.

Furthermore, the proposed exact optimization is beneficial in the light of its derivative-free feature. Such a property subject to special cases has a long history in literature while attracts renewed interest in recent years, for instance,
\cite{DuJoWaWi:15} that highlights its necessity in the very beginning. In other words, the feature is a merit since it does not require any gradients. More specifically, it is desirable in settings in which exact gradient/first-order calculations are computationally expensive or impossible. Note that, in such cases, the classical techniques -- including Kiefer-Wolfowitz-type procedures \cite{KuYi:03} -- often resort to compromises that are based on approximations, which (observation) further strengthens the above mentioned advantage in accuracy while reminds a special consideration in Remark \ref{Rem_QP_Singleton_F}. Moreover, in the fields of statistics and machine learning \cite[Refs. 3-7]{DuJoWaWi:15}, there has been a growing spectrum of optimization problems that entail functional information only, which also broadens the impacts of this research direction. Finally, it is worth reiterating that the proposed scheme exploits all the information of the function, which (unified framework) includes both the zero-order \cite{DuJoWaWi:15} and first-order values that have been widely leveraged in the classical/iterative optimization \cite{NoWr:06,Lu(Ye):03(16)}.

All of the above advantages will be validated in Sec. \ref{Sec_Ex} below, which provides explicit computational evidence. Promisingly, there exist more advantages, such as the one considered in Remark \ref{Rem_QP_Singleton_F} (Appendix \ref{App_Proof_QP}); and those of practical importance: computational speed, hardware complexity, and safety standard, since this research direction toward an exact solution aligns with, for example, \cite{JoFoTo:05} and its mathematical preliminaries (notably \cite{BeMoDuPi:02}, which is with respect to a positive definite Hessian matrix). Last but not least, an additional impact to the scientific literature suggests the (bi)linear matrix inequalities-based design (to name one recent representative, \cite{MiHa:acc}), which traces its core value/concept to the QP problem and thus motivates an interesting further research.

\section{Three Illustrative Examples}
\label{Sec_Ex}
\begin{figure}[htbp]
    \begin{center}
    \includegraphics[width=6.5cm]{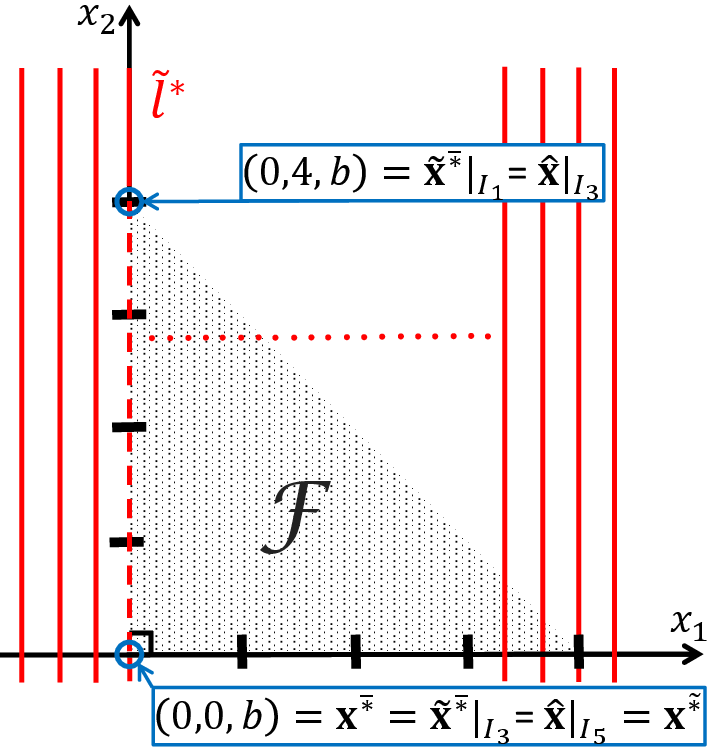}
    \caption{A geometric interpretation of ``terminal optima'' in Example {\rm\ref{Ex_QP_Singular}}, where $x_3=b$.}
    \label{Fig_Ex_QP_Singular}
    \end{center}
\end{figure}

The main objective is to demonstrate the \textit{autonomy}, \textit{generality}, and \textit{functionality} of the proposed novel \textit{closed-form} QP solvers (Theorem \ref{Thm_Equality_Constrained_QP}, Algorithm \ref{Alg_QP}, and so on), particularly from the implementation perspective. More specifically, it is worth reiterating that all the proposed results do not require any knowledge of a feasible point -- at any time before and during the solving process (namely, autonomy); deal with the general QP problem in accordance with, for example, \cite{Lu(Ye):03(16)}, as well as an extension that enlarges the application coverage (generality); and exploit all the information of the function that includes the commonly utilized zero-order \cite{DuJoWaWi:15} and first-order functional values in the classical/numerical optimization \cite{NoWr:06} (functionality). Example \ref{Ex_QP_Nonsingular} is directly adopted from \cite{Lu(Ye):03(16)}, and used to demonstrate Theorems \ref{Thm_Unconstrained_QP}, \ref{Thm_Equality_Constrained_QP}, Algorithm \ref{Alg_QP}, Corollary \ref{Cor_Alg_QP}, Remarks \ref{Rem_Unconstrained_QP}, \ref{Rem_Thm_Equality_Constrained_QP}, \ref{Rem_Cor_Alg}, \ref{Rem_Candidates_Cardinality}, and Lemma \ref{Lem_Solutions}. One step forward, Example \ref{Ex_QP_Singular} is extended and modified from Example \ref{Ex_QP_Nonsingular}, by adding an equality constraint, rendering the Hessian matrix singular, and increasing the system dimension. Therefore, it also demonstrates those mentioned above, except the special cases of Corollary \ref{Cor_Alg_QP} and Remark \ref{Rem_Cor_Alg}, and additionally Remark \ref{Rem_Terminal_Optima} for the concept of \textit{terminal optima} and related solution freedom. For compactness, the illustration for Example \ref{Ex_QP_Nonsingular} can be either inferred from Fig. \ref{Fig_Candidates} or referenced to \cite[2nd edn.,~Fig.~14.1]{Lu(Ye):03(16)}; whereas, in Example \ref{Ex_QP_Singular}, a geometric interpretation for Remark \ref{Rem_Terminal_Optima} is given in Fig. \ref{Fig_Ex_QP_Singular}. Last but not least, Example \ref{Ex_MATLAB} leverages a benchmark problem as selected in \cite{MATLAB} and compares the proposed approach with the built-in QP solver ``quadprog'', which (among others) boasts a wide popularity in the optimization literature and beyond \cite{TsMa:21}. Note that all computations in the three examples are performed and/or verified on the MATLAB$\textsuperscript{\textregistered}$ platform, and they are subject to a variety of hardware implementations (details in Table \ref{Table_Ex_MATLAB} later).

\begin{Example}
Consider a QP problem of nonsingular Hessian matrix and under inequality constraints only, where $\itn=2$, $P(1,:)=[4,1]$, $P(2,:)=[1,2]$, $\bfq=[-12,-10]^T$, $s=0$, $\bfc_1=[1,1]^T$, $\bfc_2=[-1,0]^T$, $\bfc_3=[0,-1]^T$, $d_1=4$, and $d_2=d_3=0$ \cite{Lu(Ye):03(16)}. Given $P\succ 0$, according to Corollary \ref{Cor_Alg_QP}, \ref{Thm_Unconstrained_QP_finite_value}) of Theorem \ref{Thm_Unconstrained_QP}, or line \ref{Alg_a_optimum_particular} of Algorithm \ref{Alg_QP}, we have the unique optimum with respect to the unconstrained problem as $\bfx^*=-P^{-1}\bfq=[2,4]^T$, which does not satisfy the inequality constraint associated with $(\bfc_1,d_1)$. That is, this unconstrained optimum does not reside in the feasible set. Therefore, following Algorithm \ref{Alg_QP}, the next step is the for-all environment starting at line \ref{Alg_For_All}, where $\calI=\{1,2,3\}$ ($\kappa=3$). Note that the significant parameters and values in this step are summarized by Table \ref{Table_Ex_QP_Nonsingular} in Appendix \ref{App_Tables} for completeness and easy comparison. Moreover, the examination order for all $\calI_j$'s is arbitrary, as long as every nonempty subset of $\calI$ is examined, which is discussed in \ref{Rem_Alg_Order_For_All}) of Remark \ref{Rem_Alg}.

In Table \ref{Table_Ex_QP_Nonsingular} (Appendix \ref{App_Tables}), \textit{at first} we note that only those $\calI_j$'s with $\tilde\bfb\in\calR(\tilde A)$ in the corresponding augmented system are included, as filtered at line \ref{Alg_bcd_Condition_Consistency} in Algorithm \ref{Alg_QP}; otherwise, it does not associate with the optimality since being inconsistent, and the only case is $\calI_7$, where, equivalently,
\beq
\mbox{det}(\left[\def\arraystretch{1.2} \begin{array}{c:c}
\begin{array}{c} \bfc_1^T~\\ \bfc_2^T~\\ \bfc_3^T~ \end{array} &
\begin{array}{c} ~d_1\\ ~d_2 \\ ~d_3 \end{array} \end{array}\right])\ne 0.\nonumber
\eeq

\textit{Second}, regarding $\calI_1$, the corresponding augmented matrix $\tilde A$ is of full rank, and the augmented system $\tilde A\bfx=\tilde b$ is underdetermined and has more than one solution. This falls into the processes at lines \ref{Alg_bcd_A_V_2}-\ref{Alg_bcd_Condition_Full_Rank_Singular_End} of Algorithm \ref{Alg_QP}. According to Theorem \ref{Thm_Equality_Constrained_QP}, an example of $\tilde V_2=[1/\sqrt{2},-1/\sqrt{2}]^T$ is by means of SVD, which is computed at line \ref{Alg_bcd_A_V_2} of Algorithm \ref{Alg_QP}. This yields the value ``$\tilde V_2P\tilde V_2=2\ne 0$'', which implies that this case is equivalently an unconstrained QP. As a matter of fact, given $P\succ 0$, the following results also endorse Corollary \ref{Cor_Alg_QP}:
\begin{enumerate}
\item[$\bullet$] $\tilde V_2^T(\bfq+P\tilde A^\dagger\tilde b)\in\calR(\tilde V_2 P\tilde V_2)=\calR(2)=\real$, which means that the equivalent QP problem is of finite value by {\ref{Thm_Equality_Constrained_QP_Unconstrained_QP_Finite_Condition})} of Theorem \ref{Thm_Equality_Constrained_QP}.
\item[$\bullet$] $\calN(P)=\{\bfzero\}$, and thus $\calN(\tilde A)\cap\calN(P)=\{\bfzero\}$. According to {\ref{Thm_Equality_Constrained_QP_Unconstrained_QP_Uniqueness})} of Theorem \ref{Thm_Equality_Constrained_QP} or \ref{Rem_Cor_Alg_P_Nonsingular}) of Remark \ref{Rem_Cor_Alg}, we have that the associated optimum for this case $\calI_1$ is unique and denoted $\tilde\bfx^{\bar *}=\tilde\bfx_p^{\bar *}=[1.5,2.5]^T$. This point is computed at line \ref{Alg_Inequality Equality_QP_Optimum_Condition_Inequality} in Algorithm \ref{Alg_QP}.
\item[$\bullet$] $\tilde\bfx^{\bar *}$ is an optimality candidate for the original problem, as recognized at line \ref{Alg_bcd_Condition_Full_Rank_Singular_Inclusion}, since it satisfies the other (inequality) constraints with respect to $\{(\bfc_i,d_i)\}_{i=2,3}$. In other words, it is (also) a feasible point with regard to the original problem/constraints.
\item[$\bullet$] The associated optimal value ($\tilde{\bar l}^*=-28.5$) for this case $\calI_1$ is determined by Eq. (\ref{Optimal_Value_Equality_Constrained_QP_Unconstrained_QP}), computed at line \ref{Alg_bcd_Condition_Full_Rank_Singular_Optimal_Value}, and included at line \ref{Alg_bcd_Condition_Full_Column_Rank_Inclusion} for the overall comparison later at line \ref{Alg_bcd_Result} of Algorithm \ref{Alg_QP}.
\end{enumerate}
Similar considerations apply to Cases $\{\calI_j\}_{j=2,3}$, but note that the associated unique optimum $\left(\tilde\bfx^{\bar *}=[0,5]^T\right)$ for Case $\calI_2$ is not a candidate, because it violates the first inequality constraint ``$x_1+x_2\le 4$''.

\textit{Third}, Case $\calI_4$ corresponds to a ``vertex'', which is zero-dimensional, as a result of rank$(\tilde A)=\itn=2$. Specifically, according to \ref{Rem_Thm_Equality_Constrained_QP_Vertex}) of Remark \ref{Rem_Thm_Equality_Constrained_QP}, this case is \textit{instead} included at lines \ref{Alg_bcd_Condition_Full_Column_Rank}-\ref{Alg_bcd_Condition_Full_Rank_End} of Algorithm \ref{Alg_QP}, and there is only one feasible point $\hat\bfx=[0,4]^T$ as computed at line \ref{Alg_bcd_Condition_Full_Column_Rank_Check_x}. Since $\hat\bfx$ satisfies all the constraints, it is recognized as a candidate (line \ref{Alg_bcd_Condition_Full_Column_Rank_Inclusion}). The associated, only level set value $\hat l=-24$ is exclusively optimal in this case, computed at line \ref{Alg_bcd_Condition_Full_Column_Rank_Optimal_Value}, and thus included for the comparison at line \ref{Alg_bcd_Result}. Similar considerations apply to Cases $\{\calI_j\}_{j=5,6}$, but omitted for conciseness.\vspace{0.1cm}

\textit{Finally}, all the values, either $\tilde{\bar l}^*$ or $\hat l$, are compared at line \ref{Alg_bcd_Result}. This concludes that Case $\calI_1$ corresponds to the optimality. To be more detailed, $\bfx^{\tilde *}=\tilde\bfx^{\bar *}=[1.5,2.5]^T$ for $\calI_1$ is also the unique optimum for the original problem, with the optimal value $\tilde l^*=\tilde{\bar l}^*=-28.5$, and $\bfx^{\tilde *}$ is located on the \textit{edge} of $\calF$ that is associated with $\calI_1$, explicitly, $x_1+x_2=4$.

\begin{Remark}
The results agree with \cite{Lu(Ye):03(16)} that uses ASM, but without any prior knowledge of a feasible point to initiate the solving process. Besides, all computations are performed in terms of analytical, algebraic, closed-form representations. Moreover, there are only seven cases considered in the for-all environment of Algorithm \ref{Alg_QP} (lines \ref{Alg_For_All}-\ref{Alg_For_All_End}). In consistency with Remark \ref{Rem_Candidates_Cardinality}, the cardinality of $\bar\calL$, which equals 5, is less than the upper bound ``$2^\kappa=8$''.
\label{Rem_Ex_QP_Nonsingular}
\end{Remark}

\label{Ex_QP_Nonsingular}
\end{Example}

\begin{Example}
This example is extended and modified from Example \ref{Ex_QP_Nonsingular}. The main difference in the motivation is additionally to demonstrate the concept of \textit{terminal optima}, as explained in Appendix \ref{App_Proof_QP_bcd}, and a potential to assist the analytical and comprehensive representation of all the nonunique optima. According to the problem formulation in QP (\ref{Problem_QP}), we let $\itn=3$, \textit{singular} $P=\mbox{diag}(1,0,0)\succeq 0$, $\bfq=\bfzero$, $A=[0~0~1]$, $\bfb=b\in\real$, while the other parameters follow Example \ref{Ex_QP_Nonsingular}. For succinctness, only the significantly different analyses and results are presented below.

\textit{First} (lines \ref{Alg_A_V_2}-\ref{Alg_a_QP_Condition_End} in Algorithm \ref{Alg_QP}), consider the QP under only the equality constraint ($A\bfx=b$). At line \ref{Alg_A_V_2}, $V_2$ is computed as $V_2(:,1)=[0~1~0]^T$ and $V_2(:,2)=[1~0~0]^T$. The associated matrix $V_2^TPV_2=\mbox{diag}(0,1)\ne 0$, which implies that this is equivalently an unconstrained QP by \ref{Thm_Equality_Constrained_QP_Condition_Unconstrained_QP}) of Theorem \ref{Thm_Equality_Constrained_QP}. In addition, given the singular $V_2^TPV_2\succeq 0$ as well as $V_2^T(\bfq+PA^\dagger\bfb)=\bfzero\in\calR(V_2^TPV_2)$, by {\ref{Thm_Equality_Constrained_QP_Unconstrained_QP_Finite_Condition})} of Theorem \ref{Thm_Equality_Constrained_QP}, the optimal value (that is, minimum) is finite. Moreover, since $\calN(A)\cap\calN(P)=\calR([0~1~0]^T)\ne\{\bfzero\}$, applying {\ref{Thm_Equality_Constrained_QP_Unconstrained_QP_Uniqueness})} of Theorem \ref{Thm_Equality_Constrained_QP} yields that the optima are nonunique. Preliminarily at this stage, Algorithm \ref{Alg_QP} only examines the particular solution of the optima ($\bfx_p^{\bar *}$, as computed at line \ref{Alg_a_optimum_particular}), while integrates the examination on its solution freedom ($V_2\bfvarepsilon^{\bar *}$) into the for-all environment later (lines \ref{Alg_For_All}-\ref{Alg_For_All_End}). The result is $\bfx_p^{\bar *}=[0,0,b]^T$, which resides in the feasible set $\calF$ of the original problem (grey/shaded area in Fig. \ref{Fig_Ex_QP_Singular}). Hence it qualifies as an optimality candidate (included at line \ref{Alg_Switch_Element}), associated with $\bar l^{*}=0$ (computed at line \ref{Alg_a_optimal_value}).

\textit{Second}, the process goes to the for-all environment in Algorithm \ref{Alg_QP}. Following the analyses in Example \ref{Ex_QP_Nonsingular}, the results are briefed and summarized by Table \ref{Table_Ex_QP_Singular} in Appendix \ref{App_Tables} for better presentation clarity. Note that:
\begin{enumerate}
\item[$\bullet$] Considering $\calI_1$, this corresponds to the \textit{terminal optimum} associated with the above mentioned solution freedom of $\bfx^{\bar *}$ (namely, $V_2\bfvarepsilon^{\bar *}$) as a result of the additional constraint with respect to $(\bfc_1,d_1)$. This can also be inferred from Fig. \ref{Fig_Ex_QP_Singular}, where the red/solid/vertical lines represent the level sets/values of the involved CQF; while the $x_2$-axis is just the linear variety $\bfx^{\bar *}$, and the terminal optima $(0,4,b)$ is owing to the additional constraint $(\bfc_1,d_1)$. This demonstration agrees with the discussions in Appendix \ref{App_Proof_QP_a}.
\item[$\bullet$] Regarding $\calI_2$, this is equivalently an optimization problem but of constant value, according to \ref{Thm_Equality_Constrained_QP_Condition_Unconstrained_LP}) and \ref{Thm_Equality_Constrained_QP_Condition_Constant}) of Theorem \ref{Thm_Equality_Constrained_QP}. Specifically, it is owing to $\tilde V_2^TP\tilde V_2=0$, $\tilde V_2^T(\bfq+P\tilde A^\dagger\tilde \bfb)=0$, and the constant value equals $s+[\bfq^T+\tilde\bfb^T(\tilde A^\dagger)^TP/2]\tilde A^\dagger\tilde\bfb=0$. Additionally, as addressed in {\ref{App_Proof_QP_bcd_Constant})} of Appendix \ref{App_Proof_QP_bcd}, the algorithm is designed to economize this examination by focusing on the associated terminal optima. Specifically, it is instead included in Cases $\{\calI_j\}_{j=4,5}$, and associated with additional constraints $\{(\bfc_i,d_i)\}_{i=1,3}$, respectively. Notably, the associated optimal value remains. Finally, a geometric interpretation is also available in Fig. \ref{Fig_Ex_QP_Singular}, where, on the $x_2$-axis (that is, the border of $\bfc_2^T\bfx\le d_2$), the equivalent problem/function bears the constant value. In its domain, the upper (resp., lower) terminal is due to the additional constraint ``$(\bfc_1,d_1)$'' (resp., ``$(\bfc_3,d_3)$''), and the discussion/categorization goes to the case $\calI_4$ (resp., $\calI_5$).
\item[$\bullet$] Noting the additional column of $\calN(P)\cap\calN(\tilde A)$ in Table \ref{Table_Ex_QP_Singular}, which is omitted in Table \ref{Table_Ex_QP_Nonsingular} (Appendix \ref{App_Tables}) because of the more convenient condition ``$P\succ 0$'' instead/previously in Example \ref{Ex_QP_Nonsingular}, the case $\calI_1$ (resp., $\calI_3$) renders this intersection only at the origin, and thus the corresponding optimum is unique, $\tilde\bfx_p^{\bar *}=\tilde\bfx^{\bar *}$, according to {\ref{Thm_Equality_Constrained_QP_Unconstrained_QP_Uniqueness})} of Theorem \ref{Thm_Equality_Constrained_QP}. This can also be inferred from Fig. \ref{Fig_Ex_QP_Singular}, where there is only one intersection between the optimal (level set) value, denoted $\tilde l^*$, and the border of $\bfc_1^T\bfx\le d_1$ (resp., $\bfc_3^T\bfx\le d_3$).
\end{enumerate}

\textit{Finally}, corresponding to the optimal value $\tilde l^*$, at line \ref{Alg_bcd_Result} in this algorithm, it concludes with the five optimality pairs/points residing in $\bar\calL\times\bar\calX\times\bar\calI$: $([0,0,b]^T,0,\emptyset)$, $([0,4,b]^T,0,\calI_1)$, $([0,0,b]^T,0,\calI_3)$, $([0,4,b]^T,0,\calI_4)$, and $([0,0,b]^T,0,\calI_5)$. That is, Algorithm \ref{Alg_QP} gives the optimal value, $\tilde l^*=0$, with respect to the original problem. All the \textit{terminal optima} are $[0,0,b]^T$ and $[0,4,b]^T$, as illustrated in Fig. \ref{Fig_Ex_QP_Singular}. Notably, the former i) is the particular solution with regard to the QP under only the equality constraint, and also feasible to the original QP; and ii) resides on an edge associated with the inequality constraint ``$\bfc_3^T\bfx\le d_3$'', while also on a vertex with the additional constraint ``$\bfc_2^T\bfx\le d_2$''. On the other hand, the latter terminal optimum resides on an edge associated with the constraint ``$\bfc_1^T\bfx\le d_1$'', while also on a vertex with the additional ``$\bfc_2^T\bfx\le d_2$''. It is worth remarking that Algorithm \ref{Alg_QP} ably gives all/both the \textit{terminal optima}. Based on this, it can be analyzed that all the \textit{optima} are just \textit{intermediate}, denoted as $\calX^*=\{[0,t,b]^T~\vrule~0\le t\le 4,~t\in\real\}$, which can also be revealed from Fig. \ref{Fig_Ex_QP_Singular} (specifically, the red-dashed line segment).

\begin{Remark}
As expected and similar to Example \ref{Ex_QP_Nonsingular} (Remark \ref{Rem_Ex_QP_Nonsingular}), Algorithm \ref{Alg_QP} autonomously yields the results, with only functional information (derivative-free) and without any (prior) knowledge of a feasible (starting) point. Also, in this setting, all computations are performed in terms of algebraic and explicit representations. Moreover, there are still only few cases (seven, to be exact) considered in the for-all environment of Algorithm \ref{Alg_QP} and, being consistent with Remark \ref{Rem_Candidates_Cardinality}, the cardinality of $\bar\calL$ (6) is indeed less than the upper bound ``$2^\kappa=8$''.
\end{Remark}

\label{Ex_QP_Singular}
\end{Example}

\begin{Example}
The extensive QP literature highlights a universal and popular solver ``quadprog'' in MATLAB$\textsuperscript{\textregistered}$ \cite{TsMa:21}, and this example leverages a benchmark problem as demonstrated in \cite{MATLAB} to offer more computational evidences of the superiority of the proposed approach. Specifically, let $P(1,:)=[1,-1,1]$, $P(2,:)=[-1,2,-2]$, $P(3,:)=[1,-2,4]$, $\bfq=[-7,-12,-15]^T$, $s=0$, $A=[1,1,1]$, and $b=3$ in Problem/Eq. (\ref{Problem_Equality_QP}), and this Example \ref{Ex_MATLAB} compares Theorem \ref{Thm_Equality_Constrained_QP} with ``quadprog''. Among various implementations for Theorem \ref{Thm_Equality_Constrained_QP}, we demonstrate an intuitive and preliminary one using MATLAB$\textsuperscript{\textregistered}$ for the sake of brevity, whose primary steps are sequentially summarized:

\vspace{0.6cm}\hspace{0.3cm}\begin{minipage}{0.95\textwidth}
\begin{lstlisting}[frame=trBL, numbers=left, escapeinside={(*@}{@*)}, escapechar=|]
A_dagger=pinv(A);|\label{Ex_MATLAB_Implementation_pinv}|
V_2=null(A);|\label{Ex_MATLAB_Implementation_null}|
x_star_particular=A_dagger*b...|\label{Ex_MATLAB_Implementation_Optimum_1}|
  -V_2*pinv(V_2'*P*V_2)*(V_2'*q+V_2'*P*A_dagger*b);|\label{Ex_MATLAB_Implementation_Optimum_2}|
l_star_equilibrium=(b'*A_dagger'*P/2+q')*(eye(n)-V_2*pinv(V_2'*P*V_2)...|\label{Ex_MATLAB_Implementation_Optimal_Value_1}|
  *V_2'*P)*A_dagger*b+s-q'*V_2*pinv(V_2'*P*V_2)*V_2'*q/2;              (i)|\label{Ex_MATLAB_Implementation_Optimal_Value_2}|
\end{lstlisting}
\end{minipage}
\vspace{0.2cm}

\noindent Note that lines \ref{Ex_MATLAB_Implementation_Optimum_1} and \ref{Ex_MATLAB_Implementation_Optimum_2} in Implementation (i) correspond to Eq. (\ref{Optimum_Equality_Constrained_QP_Unconstrained_QP}) while \ref{Ex_MATLAB_Implementation_Optimal_Value_1} and \ref{Ex_MATLAB_Implementation_Optimal_Value_2} to (\ref{Optimal_Value_Equality_Constrained_QP_Unconstrained_QP}). In addition, this example falls into Case \ref{Thm_Equality_Constrained_QP_Condition_Unconstrained_QP}) of Theorem \ref{Thm_Equality_Constrained_QP} because of $\calR(P)\cap\calN(A)=\calN(A)\ne\{\bfzero\}$, which can be easily verified -- using $V_2^TPV_2\ne 0$ -- and thus omitted in the summary of primary steps (i).

On the MATLAB$\textsuperscript{\textregistered}$ platform, the proposed solver yields \textit{the} optimal value $\bar l^*=-47.1786$ and \textit{the} optimum $\bfx^{\bar *}=[-3.5714,2.9286,3.6429]^T$, where the uniqueness of $\bar l^*$ (resp., $\bfx^{\bar *}$) is owing to the nonsingular $V_2^TPV_2$, namely the nulled flexibility $\bfvarepsilon^{\bar *}=\bfzero$ in Eq. (\ref{Optimum_Equality_Constrained_QP_Unconstrained_QP}); whereas, the built-in ``quadprog'' routine renders the same $\bar l^*$ and $\bfx^{\bar *}$ but \textit{not} the associated uniqueness, respectively \cite{MATLAB}. Moreover, the performance in computation time also quite differ, which is evaluated -- for better comprehensiveness -- subject to a variety of software and hardware in Table \ref{Table_Ex_MATLAB}. Specifically, Table \ref{Table_Ex_MATLAB} also more favors Theorem \ref{Thm_Equality_Constrained_QP} than ``quadprog'', where the evaluation is in terms of the ratio of computation time (Theorem \ref{Thm_Equality_Constrained_QP} over ``quadprog'') while among several Intel Core$\textsuperscript{\textregistered}$ CPUs and MATLAB$\textsuperscript{\textregistered}$ versions till the state of the art (with 8GB RAM). Accordingly, the maximum (resp., minimum) advantage amounts to 2.3$\%$ (resp., 9.6$\%$) of the computation time using ``quadprog''; in other words, the proposed solver reduces the computational effort by 90.4$\%-$97.7$\%$. In average, the numerics still significantly account for a computational reduction of 94.5$\%$ by virtue of Theorem \ref{Thm_Equality_Constrained_QP}. It is worth reiterating that the Implementation (i) is only preliminary, whose efficiency potentials can be more exploited by, for instance, ``coalescing'' lines \ref{Ex_MATLAB_Implementation_pinv} and \ref{Ex_MATLAB_Implementation_null}, namely, by performing the SVD operation on $A$ only once. This requires modifications for the ``svd'' routine in MATLAB$\textsuperscript{\textregistered}$, and thus serves as a future research to remain focused of this article.

\begin{table}[htbp]
\centering
\begin{threeparttable}[htbp]
\caption{Advantage in computation time\tnote{a} in Example {\rm\ref{Ex_MATLAB}}}
\renewcommand{\arraystretch}{1.5}
\begin{tabular}{l c c c}
\hline
 & ~~~i3, 2.75GHz & ~~~i5, 3.4GHz & ~~i7, 2.2GHz \\ \hline
MATLAB$\textsuperscript{\textregistered}$ 2010b & ~~~2.5$\%$ & ~~~2.4$\%$ & ~~2.3$\%$ \\
MATLAB$\textsuperscript{\textregistered}$ 2020b & ~~~8.9$\%$ & ~~~7.4$\%$ & ~~9.6$\%$ \\ \hline
\end{tabular}
\begin{tablenotes}
\item [a] Ratio of time using Theorem \ref{Thm_Equality_Constrained_QP} over ``quadprog''.
\end{tablenotes}
\label{Table_Ex_MATLAB}
\end{threeparttable}
\end{table}

\begin{Remark}
Regarding the up-to-date universal solver ``quadprog'', \cite{MATLAB} highlights a \textit{necessary} condition ``$P\succ 0$'' that associates with \textit{a} finite optimal value, even if the QP problem is without any constraint (Sec. \ref{Subsec_Unconstrained_QP}). Actually, there exists a counterexample to the condition/necessity: $P=\mbox{diag}(1,0)$, $\bfq=\bfzero$, and $s=0$ in Eq. (\ref{CQF}). That is, the QP problem is to minimize $F_\bfx(\bfx)=x_1^2$, which has a/the finite minimum (zero) whereas the Hessian matrix is not positive definite. As a comparison, \ref{Thm_Unconstrained_QP_finite_condition}) of Theorem \ref{Thm_Unconstrained_QP} provides a \textit{necessary and sufficient} condition, which readily guarantees \textit{the} finite minimum in the above counterexample -- owing to $\bfq=\bfzero\in\calR(P)$. Remarkably, Theorem \ref{Thm_Unconstrained_QP} fundamentally helps develop the remaining results (with constraints) in Sec. \ref{Sec_Convex_Optimization}, such as another \textit{equivalent} condition for \textit{the} finite optimal value in {\ref{Thm_Equality_Constrained_QP_Unconstrained_QP_Finite_Condition})} of Theorem \ref{Thm_Equality_Constrained_QP}.
\end{Remark}

\label{Ex_MATLAB}
\end{Example}

\section{Conclusion}
\label{Sec_Conclusion}
A motivation of this application -- to the field of nonlinear programming based on the recent observations of CQE -- responds to an open question in \cite{DuJoWaWi:15}; and an expectation in the general optimization literature: computational enhancement for practical advances \cite{Lu(Ye):03(16)}. The philosophy aims at analytical and exact optimization, as a comparison to the common numerical algorithms \cite{VaXiYa:20,NoWr:06}. Specifically, the focus of this article is on the QP, which is a constituent in various approaches such as the Newton's method. At first, we re-examine the unconstrained QP based on, and from the perspective of, the new analyses of CQE (in the recent literature), CQF, and CQE-CQF relation. This preliminary finding facilitates a complete and analytical characterization of the equality-constrained QP, which actually can be categorized into three equivalent problems in a unified framework. In addition, all the above results are consistent with the literature, when specialized to specific considerations; while promising in terms of the computational performance, which is partially owing to the shared motivation and beneficial results in the recent contribution \cite{LiLiHs:20}. Another highlight in this presentation is the proposed QP solver, in accordance with the general problem formulation. The analytical philosophy also benefits from the above mentioned categorization of equality-constrained QP, and yields algebraic zero-order closed-form results/solutions in a guaranteed, pre-determined, finite, and explicit number of steps; and without any knowledge of a feasible point, a priori and any time in the process. All these characteristics/merits are shared throughout this article. Moreover, to further exploit its computational capabilities, this flexible QP solver includes problem-specific (and more efficient) variants below:
\begin{enumerate}
\item[$\bullet$] an economic/expedited version for a subset of QP problems;
\item[$\bullet$] a special but popular case ``the objective function whose Hessian matrix is positive definite'', which corresponds to a thinned version of the algorithm;
\item[$\bullet$] QP problems subject to inequality constraints only, being applicable (that is, solvable) after a simple replacement of parameters in the algorithm.
\end{enumerate}
Note that, regarding a subset of QP problems, we preliminarily and explicitly present the associated terminal optimum/optima. This indicates a further research direction, because it requires an analytical representation of all the elements bounded by the inequality constraints. Moreover, we extend the general formulation on the constraints to a different branch, as inspired by the analyses of CQE and CQF, which further broadens the spectrum of QP applications. Last but not least, two representative solvers in literature (the ASM and barrier method) help justify advantages/potentials of the proposed results; while three examples demonstrate the results with illustrations, computational evidences, and further insights (notably, a \textit{counterexample} to the state of the art) -- analytically as well as numerically on the MATLAB$\textsuperscript{\textregistered}$ platform. Under positive definite and semidefinite Hessian matrices in the objective functions, respectively, these demonstrations endorse the effectiveness, autonomy, efficiency, functionality, and exactness of the proposed solvers.

\appendices
\section{Proof of Theorem \ref{Thm_Unconstrained_QP} (Unconstrained QP)}
\label{App_Proof_Unconstrained_QP}

For brevity, only the results that are significantly required in this proof are given, whereas most of the shared similarities are referenced to Lemma \ref{Lem_Solutions} and its proof.

\noindent\ref{Thm_Unconstrained_QP_all_preimage}) To obtain the preimage of any level set, $F_\bfx=l$, is equivalent to that of the zero level set of $\tilde F_\bfx$, where the CQF $\tilde F_\bfx:\real^\itn\rightarrow\real$ and $\tilde F_\bfx=F_\bfx-l$, as ``(vertically) shifted'' from $F_\bfx$ by $l$. Furthermore, this is equivalent to solving the CQE, $\bfx^TP\bfx+\bfq^T\bfx+s-l=0$, which is always solvable since the original process is to find the preimage. The result directly follows by applying Lemma \ref{Lem_Solutions} to this CQE. Note that the above mentioned ``shift'' is widely used throughout Sec. \ref{Sec_Convex_Optimization}, and thus exemplified/illustrated in Fig. \ref{Fig_Shift} for presentation clarity.\vspace{0.16cm}

\noindent\ref{Thm_Unconstrained_QP_finite_condition}) Divide the derivations, following Lemma \ref{Lem_Solutions}, into whether $P$ is of full rank. A) If rank$(P)=\itn$, then we reformulate the CQF $F_\bfx$ in Eq. (\ref{CQF}) equivalently as
\beq
F_\bfx(\bfx)=\left\Vert P^{1/2}\bfx+P^{-{1/2}}\bfq\right\Vert^2/2+s-\bfq^TP^{-1}\bfq/2,
\label{Proof_Thm_Unconstrained_QP_CQF_n}
\eeq
which is similar to \cite[Eq.~(34)]{LiLiHs:20}. Obviously, this function (\ref{Proof_Thm_Unconstrained_QP_CQF_n}) is bounded below by the value ``$(s-\bfq^TP^{-1}\bfq/2)$'', which (value) is finite and important in the derivations afterward. Note that, in this case, $\bfq\in \calR(P)=\real^\itn$. On the other hand, B) if rank$(P)=r<\itn$, then reformulate $F_\bfx$ (\ref{CQF}) equivalently as
\beq
F_\bfx(\bfx)=\bfx_1^T\breve\Sigma_1\bfx_1+\bfq_1^T\bfx_1+s+\bfq_2^T\bfx_2,
\label{Proof_Thm_Unconstrained_QP_CQF_r}
\eeq
which is similar to \cite[Eq.~(42)]{LiLiHs:20}. Obviously, if $\bfq_2=\bfzero$, that is, $\bfq\in\calR(P)$, then the domain of $F_\bfx$ \textit{shrinks} and corresponds to the $\bfx_1$-freedom in $\bfx$. Similar to A) above, $F_\bfx$, or equivalently the CQF $F_{\bfx_1}:\real^\itr\rightarrow\real$,
\beq
F_{\bfx_1}(\bfx_1)=\bfx_1^T\breve\Sigma_1\bfx_1+\bfq_1^T\bfx_1+s,
\label{Proof_Thm_Unconstrained_QP_CQF_r_x1}
\eeq
is bounded below by $\breve l\coloneqq s-\bfq_1^T\breve\Sigma_1\bfq_1/2=s-\bfq^T P^\dagger\bfq/2\in\real$. However, if $\bfq_2\ne\bfzero\Leftrightarrow \bfq\not\in\calR(P)$, then $F_\bfx$ in Eq. (\ref{Proof_Thm_Unconstrained_QP_CQF_r}) is unbounded because of the $\bfx_2$-freedom, which is decoupled from $\bfx$.\vspace{0.16cm}

\noindent\ref{Thm_Unconstrained_QP_finite_value}) Given $\bfq\in \calR(P)$, the derivations differ only in the singularity of $P$ according to Conditions/Eqs. (\ref{Solvability_Cond_n}) and (\ref{Solvability_Cond_k_R(M)}), which are divided into Cases A) and B) below.\vspace{0.16cm}

\noindent A) If rank$(P)=\itn$, by \ref{Thm_Unconstrained_QP_all_preimage}) above and Eq. (\ref{Solution_n}), the unique optimum happens exclusively when the only solution freedom ($\bfv$) is canceled. In other words, the term in the square-root operator (``$-l^*$'' as in this theorem) vanishes. Geometrically from Fig. \ref{Fig_Shift}, this can also be explained by the \textit{critical shift} of $l^*$, which allows the only one-point intersection with the zero plane/line. It is worth mentioning that this shift can be ``ascending'' (resp., ``descending''), if $l^*<0$ (resp., $l^*>0$), that is, the critical ``distance'' to \textit{null} the square-root operation (resp., make it \textit{consistent}). For succinctness, another explanation using Eq. (\ref{Proof_Thm_Unconstrained_QP_CQF_n}) is omitted. To summarize, the preimage of the level set at the optimal value $l^*$ is the singleton $\{\bfx^*\}$.\vspace{0.16cm}

\noindent B) If rank$(P)=r<\itn$, we at first introduce the concept of \textit{hierarchical layers} in the preimage of CQF in Eq. (\ref{CQF}), according to \ref{Thm_Unconstrained_QP_all_preimage}) above and Eq. (\ref{Solution_k_R(M)}). There are three layers: I) the unique point, $\bfx_p^*=-P^\dagger\bfq$, II) the $\bfvarepsilon$-freedom in $\calN(P)$, and III) the $\bfrho$-freedom in $\calR(P)$. This I)-III) order follows the size of included preimage elements, and I) being the smallest. Geometrically, the I) and II) layers correspond to the optima, where II) mainly reflects the singularity of $P$. Then, introducing III) equips with the full freedom, which thus corresponds to the complete preimage. For brevity, other explanations by means of, for instance, a) the geometrical perspective, which interprets the critical shift of optimal value $l^*$ similarly as in A); b) Eq. (\ref{Proof_Thm_Unconstrained_QP_CQF_r_x1}), which can be reformulated similar to (\ref{Proof_Thm_Unconstrained_QP_CQF_n}); and c) the support from \cite[Theorem~ 3.2]{LiLiHs:20} (and its proof) are omitted.

\begin{Remark}
As an example, Fig. \ref{Fig_Layers} illustrates the concept of ``hierarchical layers''. Regarding the singular $P$, the three layers are represented by the I) black/crossed point: this particular solution is denoted by $\bfx_p^*=-P^{-1}\bfq$ at the center of the figure, II) blue/gridded closed region, and III) red/dashed region, respectively. All the layers are within the preimage of $F_\bfx$, denoted by $F_\bfx^{-1}\coloneqq F_\bfx^{-1}(l):\real\rightarrow\real^\itn$ with respect to any level set value $l$. I) and II) layers are in the same level set at (the optimal value) $l^*$, while each red/dashed ellipse in the III) layer corresponds to a level set at a higher value. Moreover, II) corresponds to $\calN(P)$ and each element/point is referenced by the vector $\breve\bfepsilon$; while III) to a point $\Phi\cdot\breve\bfrho\in\calR(P)$, where $\Phi\coloneqq\Phi(l)=\sqrt{2l+\bfq^TP^\dagger\bfq-2s}\cdot P^{\dagger/2}\in\real^{\itn\times\itn}$ while $\breve\bfrho\in\calR(P)\subseteq\real^\itn$ is of unit length. Note that $\calR(\Phi)=\calR(P^{\dagger/2})=\calR(P)$ and $\Vert\Phi\Vert=\sqrt{(2l+\bfq^TP^\dagger\bfq-2s)/\sigma_\itr}$ \cite{GoVa:13}, where $\sigma_\itr>0$ denotes the smallest (positive) singular value of $P$. On the other hand, considering the case of invertible $P$, the only differences are
\begin{enumerate}
\item[$\bullet$] the lack of II) layer, because of $\bfvarepsilon\in\calN(P)=\{\bfzero\}$ and $P^\dagger=P^{-1}$;
\item[$\bullet$] $\Phi=\sqrt{2l+\bfq^TP^{-1}\bfq-2s}\cdot P^{-1/2}$ and $\Vert\Phi\Vert=\sqrt{(2l+\bfq^TP^{-1}\bfq-2s)/\sigma_\itn}$, where $\sigma_\itn$ denotes the smallest (positive) singular value or, equivalently here, eigenvalue of $P$.
\end{enumerate}
Finally, this concept of ``hierarchical layers'' is essential for the subsequent derivations and results in Sec. \ref{Sec_Convex_Optimization}, particularly the constrained optimization as analytically solved from a novel perspective (Secs. \ref{Subsec_Equality_QP} and \ref{Subsec_QP}).
\label{Rem_Fig_Geometric_Layers}
\end{Remark}

\section{Proof of Theorem \ref{Thm_Equality_Constrained_QP} (Equality Constraints)}
\label{App_Proof_Equality_Constrained_QP}

Let the SVD of $A$ be given by
\beq
A=W\cdot\left[ \begin{array}{cc} \Gamma_1~~ & 0 \end{array}\right]\left[ \def\arraystretch{1.2}\begin{array}{c} V_1^T\\V_2^T\end{array}\right],\nonumber
\eeq
with its unique pseudoinverse $A^\dagger\in\real^{\itn\times\itm}$ \cite{GoVa:13},
\beq
A^\dagger&=&[V_1~V_2]\left[ \begin{array}{c} \Gamma_1^{-1}\\ 0 \end{array}\right]\cdot W^T\nonumber\\
&=&V_1\Gamma_1^{-1}W^T~\mbox{(the thin version)},
\label{App_Proof_Equality_Constrained_QP_A_dagger}
\eeq
where $W\in\real^{\itm\times\itm}$, $V_1\in\real^{\itn\times\itm}$, $V_2\in\real^{\itn\times(\itn-\itm)}$, $\calR(V_1)=\calR(A^T)$, $\calR(V_2)=\calN(A)$, $\Gamma_1\in\real^{\itm\times\itm}$ is nonsingular, and more properties are widely available in literature \cite{GoVa:13}. The idea is to equivalently transform the equality-constrained QP into an unconstrained optimization problem, by down-casting to the linear variety/constraint as described by $A\bfx=\bfb$. Specifically, the first step is to parameterize all the feasible points, namely $A\bfx=\bfb\Leftrightarrow$ Eq. (\ref{Parameterization_Equality_Constraint}), where $\bfy\in\real^{\itn-\itm}$. Then, reformulate the considered constrained optimization problem, ``minimize $F_\bfx$ (\ref{CQF}) subject to $A\bfx=\bfb$'', equivalently as the following one: minimize $F_\bfy$ with respect to $\bfy$, where $F_\bfy:\real^{\itn-\itm}\rightarrow\real$,
\beq
F_\bfy(\bfy)=\bfy^TV_2^TPV_2\bfy/2+[\bfb^T(A^\dagger)^TP+\bfq^T]V_2\bfy+s+[\bfq^T+\bfb^T(A^\dagger)^TP/2]A^\dagger\bfb.
\label{Proof_Equality_Constrained_QP_Equivalent_y}
\eeq
Therefore, the QP Problem (\ref{Problem_Equality_QP}) is equivalently categorized into the three cases \ref{Thm_Equality_Constrained_QP_Condition_Unconstrained_QP})-\ref{Thm_Equality_Constrained_QP_Condition_Constant}), according to the coefficients of quadratic, linear, and constant terms in (\ref{Proof_Equality_Constrained_QP_Equivalent_y}), respectively. The system dimension is expectedly reduced, from $\real^\itn$ to $\real^{\itn-\itm}$. The following derivations focus on the case of (equivalently unconstrained) QP, while it is straightforward to derive the other cases and thus omitted for compactness.

To start with, we further analyze the equivalence conditions such that $F_\bfy$ in Eq. (\ref{Proof_Equality_Constrained_QP_Equivalent_y}) is quadratic, intentionally/preferably in terms of the original, given parameters. Let the SVD of (non)singular $P\succeq 0$ be
\beq
P=\check U_1\check\Sigma_1\check U_1^T~\mbox{(full/thin version)},
\label{App_Proof_Equality_Constrained_QP_P_SVD}
\eeq
where $\check U_1\in\real^{{\it n}\times {\check{\it r}}}$ is with orthonormal column(s), $\check\itr\coloneqq\mbox{rank}(P)$, and $\check\Sigma_1\in\real^{{\check\itr}\times {\check\itr}}$ is nonsingular. With respect to the case of singular $P$, more properties and details regarding this SVD similarly follow \cite[Eqs.~{(35)-(37)}]{LiLiHs:20}, which can be easily extended to the nonsingular/other case. It is worth emphasizing that the following derivations \textit{coalesce} both cases of $P$. On the one hand, this is a QP, namely, $F_\bfy$ in Eq. (\ref{Proof_Equality_Constrained_QP_Equivalent_y}) is a CQF, iff $V_2^TPV_2\ne 0\Leftrightarrow V_2^T\check U_1\check \Sigma_1\check U_1^TV_2\ne 0\Leftrightarrow \check U_1^T V_2\ne 0\Leftrightarrow \calR(P)\cap\calN(A)\ne\{\bfzero\}$. Similarly on the other, this is \textit{not} a QP, iff $\check U_1^T V_2=0\Leftrightarrow \mbox{``}\calN(A)\subseteq(\calR(P))^\perp=\calN(P)$'', where the equality holds when $\calR(P)\oplus\calN(A)=\real^\itn$.

Notably, to conform to the applicable form using Theorem \ref{Thm_Unconstrained_QP}, any level set of $F_\bfy$ (\ref{Proof_Equality_Constrained_QP_Equivalent_y}) is regarded as a CQE (\ref{CQE}) in the unknown variable $\bfz=\bfy$, with the parameters given in {\ref{Thm_Equality_Constrained_QP_Unconstrained_QP_All_Preimage})} of this theorem. Applying Theorem \ref{Thm_Unconstrained_QP} yields the results in \ref{Thm_Equality_Constrained_QP_Condition_Unconstrained_QP}), except {\ref{Thm_Equality_Constrained_QP_Unconstrained_QP_Uniqueness})}, while omitting the lengthy but somewhat straightforward calculations. Regarding the exceptional/last one {\ref{Thm_Equality_Constrained_QP_Unconstrained_QP_Uniqueness})}, based on the above derivations, the optimum is unique, iff the Hessian matrix of the CQF $F_\bfy$ in Eq. (\ref{Proof_Equality_Constrained_QP_Equivalent_y}) is positive definite. Rewrite this equivalence condition as $V_2^TPV_2=V_2^T\check U_1\check\Sigma_1\check U_1^TV_2\succ 0$, according to the thin version of Eq. (\ref{App_Proof_Equality_Constrained_QP_P_SVD}). Note that the full version corresponds to the case of nonsingular $P$, and the result similarly follows. Let any $\tilde\bfy\in\real^{\itn-\itm}\backslash\{\bfzero\}$ be given, this condition is further equivalent to
\beq
&&\Vert\check\Sigma_1^{1/2}\check U_1^TV_2\tilde\bfy\Vert^2>0\nonumber\\
&&\Leftrightarrow\check U_1^T\cdot V_2\tilde\bfy\ne\bfzero\nonumber\\
&&\Leftrightarrow\calN(P)\cap\calN(A)=\{\bfzero\},\nonumber
\eeq
where $\check\Sigma_1^{1/2}\coloneqq\mbox{diag}(\sqrt{\sigma_1},\sqrt{\sigma_2},\cdots,\sqrt{\sigma_{\check\itr}})\succ 0$ and $\{\sigma_\iti\}_{\iti=1}^{\check\itr}$ denotes the set of nonzero singular values or, equivalently in this case, eigenvalues of $P$.

\begin{Remark}
Regarding the SVD of $A$, the columns of $V_2$ are not unique. However, the effect by this nonuniqueness is completely coalesced by virtue of Eq. (\ref{Parameterization_Equality_Constraint}), in the very beginning of this proof. In other words, the results in Theorem \ref{Thm_Equality_Constrained_QP} are independent of the nonunique choices of $V_2$. From the implementation perspective, in the literature there exist many and various algorithms to construct an example, such as the Golub-Reinsch SVD till the second step \cite{GoVa:13}. Note that, in the MATLAB$\textsuperscript{\textregistered}$ platform, the ``null'' function ably gives an example of $V_2\in\calN(A)$, which is implemented by computing the SVD of $A$.
\end{Remark}

\begin{Remark}
For computational purposes, we select the SVD process to obtain a basis of $\calN(A)$, that is, the right singular vectors of $A$ that are associated with the zero singular value. As expected, the derivations only utilize the property of $\calR(V_2)=\calN(A)$ where $V_2\in\real^{\itn\times(\itn-\itm)}$, and can be easily extended to the general case: any matrix in $\real^{\itn\times(\itn-\itm)}$ with its range spanning $\calN(A)$. However, owing to many computational benefits \cite{GoVa:13} while noting the largely used parameter $A^\dagger$ (and $P^\dagger$) in Sec. \ref{Sec_Convex_Optimization}, the SVD process is more reasonably preferred to construct an example. Besides, regarding the extension in Sec. \ref{Subsec_QP}, the proposed algorithm is based on this theorem, and thus shares such a consideration.
\end{Remark}

\section{Proof of Algorithm \ref{Alg_QP} (Exact QP Solver)}
\label{App_Proof_QP}

\textit{Outline}: The proof builds on the concept ``hierarchical layers'' as illustrated in Fig. \ref{Fig_Layers}, which clarifies the geometric structure of the solution to CQE (\ref{CQE}) that is associated with the optimality/minimality of CQF (\ref{CQF}); and case-by-case arguments in accordance with Fig. \ref{Fig_Locations}, which exactly extract the optimality for most of the cases, except a subset that requires further investigation (parameterization of inequality-constrained feasible points) to complementarily locate all the optima. Moreover, the derivations largely utilize Theorem \ref{Thm_Equality_Constrained_QP}, which are divided into: Appendix \ref{App_Proof_QP_a} corresponding to lines \ref{Alg_a_QP_Condition}-\ref{Alg_a_QP_Condition_End} in Algorithm \ref{Alg_QP}; and \ref{App_Proof_QP_bcd} to lines \ref{Alg_bcd}-\ref{Alg_bcd_Result}, which augments/coalesces the inequality constraints with respect to each nonempty subset (line \ref{Alg_For_All}). Note that the design of categorizations examines every ``optimality candidate'' (a demonstration in Fig. \ref{Fig_Candidates}) for the ultimate comparison at line \ref{Alg_bcd_Result}, and benefits efficient extensions to QP variants such as those in Corollary \ref{Cor_Alg_QP}.

To start with, Fig. \ref{Fig_Locations} illustrates the main categorization for the derivations in this proof. With regard to the QP Problem (\ref{Problem_QP}) under only equality constraints, (a) considers the case that the unique optimum, or the particular solution $\bfx_p^{\bar *}$ associated with the optima, resides in the feasible set $\calF$ (shaded/grey area, a polyhedron \cite{BoVa:04,Lu(Ye):03(16)}). If the considered point in (a) is outside of $\calF$, but the degree of freedom in (\ref{Optimum_Equality_Constrained_QP_Unconstrained_QP}) associated with the case of nonunique optima ($V_2\bfvarepsilon^{\bar *}$) intersects $\calF$, then this is included in (b). On the other hand ($\bfx^{\bar *}\not\in\calF$), the derivations go to either (c) or (d), where (c) addresses the case that the optimality occurs on an ``edge/space of nonzero dimension'' while (d) at a ``vertex/point''. Note that $\calF\ne\emptyset$, according to the QP problem formulation in (\ref{Problem_QP}). Moreover, denote $\tilde\calI^*\subseteq\calI$ as the \textit{optimality subset of involved inequality constraints}, which associates with the additional equality constraint(s) that is/are the \textit{border(s)} of involved inequality constraint(s). This subset indicates the ``edge(s)'' or ``vertex/vertices'' that the optimum/optima reside(s). A special case is $\tilde\calI^*=\emptyset$, which corresponds to the case when QP Problems (\ref{Problem_Equality_QP}) and (\ref{Problem_QP}) yield the identical solution. Finally, the first line (fundamental step) of this algorithm is to compute an orthonormal basis of $\calN(A)$, denoted $V_2$, preferably using SVD due to its many computational advantages \cite{GoVa:13}.

\subsection{Case (a) $\bfx_p^{\bar *}\in\calF$}
\label{App_Proof_QP_a}
According to Theorem \ref{Thm_Equality_Constrained_QP}, this case corresponds to \ref{Thm_Equality_Constrained_QP_Condition_Unconstrained_QP}) where $V_2^TPV_2\ne 0$, as filtered at line \ref{Alg_a_QP_Condition}. Note that the case of $V_2^TPV_2=0$ corresponds to either an LP or constant function. If, due to additional inequality constraints, the optimality of QP Problem (\ref{Problem_QP}) is associated with the former (resp., latter), then the optimality occurs (resp., also) at the ``terminal optima/optimum'' that will be detailed later in Appendix \ref{App_Proof_QP_bcd} while determined/computed at lines \ref{Alg_bcd}-\ref{Alg_bcd_Result} of this algorithm.

By {\ref{Thm_Equality_Constrained_QP_Unconstrained_QP_Uniqueness})} of Theorem \ref{Thm_Equality_Constrained_QP}, the optimum of equality-constrained QP is unique, iff $\calN(A)\cap\calN(P)=\{\bfzero\}$, or equivalently $V_2^TPV_2\succ 0$ in Eq. (\ref{Optimum_Equality_Constrained_QP_Unconstrained_QP}). Accordingly, further divide this case into whether $V_2^TPV_2$ is nonsingular.\vspace{0.16cm}

\noindent i) If $V_2^TPV_2\succ 0$, then the unique optimum $\bfx^{\bar *}=\bfx_p^{\bar *}$ is given by Eq. (\ref{Optimum_Equality_Constrained_QP_Unconstrained_QP}), where the freedom $\bfvarepsilon^{\bar *}\in\calN(V_2^TPV_2)=\{\bfzero\}$ is nulled, and included in line \ref{Alg_a_optimum_particular} of this algorithm. The associated optimal value $\bar l^*$ is given by Eq. (\ref{Optimal_Value_Equality_Constrained_QP_Unconstrained_QP}), and included in line \ref{Alg_a_optimal_value}. However, only if this unique optimum resides in the feasible set $\calF$, which is verified at line \ref{Alg_a_optimum_in}, the algorithm completes solving QP Problem (\ref{Problem_QP}), since this optimality $(\bfx^{\bar *},\bar l^*)$ under only equality constraints also solves the general QP Problem (\ref{Problem_QP}). It is worth mentioning that the pair of ``$V_2^TPV_2\succ 0$ and the additional, returned parameter $\tilde\calI^*=\emptyset$'' serves as a unique identifier for this case, but this excludes the extreme case of singleton $\calF\ne\emptyset$ (details at Remark \ref{Rem_QP_Singleton_F} later).\vspace{0.16cm}

\noindent ii) If $V_2^TPV_2$ is singular, then at first we need to check whether the optimal value is finite, whose equivalence condition is formulated in {\ref{Thm_Equality_Constrained_QP_Unconstrained_QP_Finite_Condition})} of Theorem \ref{Thm_Equality_Constrained_QP} and included in line \ref{Alg_Equality_QP_Optimum_Condition_Inequality}. If infinite, then the optimality of QP Problem (\ref{Problem_QP}) can only possibly occur at the ``terminal optima/optimum'', owing to additional inequality constraints, which will be determined later in Appendix \ref{App_Proof_QP_bcd}. On the other hand, consider the finite optimal value as given by Eq. (\ref{Optimal_Value_Equality_Constrained_QP_Unconstrained_QP}), with the associated optima by (\ref{Optimum_Equality_Constrained_QP_Unconstrained_QP}). These two have also been computed in lines \ref{Alg_a_optimum_particular} and \ref{Alg_a_optimal_value} for a concise presentation. Note that, at this stage, only the particular solution of Eq. (\ref{Optimum_Equality_Constrained_QP_Unconstrained_QP}) is checked to be feasible or not; while the freedom arising from $V_2\bfvarepsilon^{\bar *}$, where $\bfvarepsilon^{\bar *}\in\calN(V_2^TPV_2)\ne\{\bfzero\}$, will be checked later in the for-all environment (at lines \ref{Alg_For_All}-\ref{Alg_For_All_End}). This argument similarly follows {\ref{App_Proof_QP_bcd_QP_bounded})} in Appendix \ref{App_Proof_QP_bcd} below. It is worth remarking that an expedited alternative is available by economizing the searching of any possible further optimum, as determined by some $\bfvarepsilon^{\bar *}\ne\{\bfzero\}$. Specifically, replacing lines \ref{Alg_Switch_Element}-\ref{Alg_Switch_Go_To} by \ref{Alg_Equality_QP_Positive_Definite_In}.

\subsection{Cases {(b)-(d)} $\bfx_p^{\bar *}\not\in\calF$}
\label{App_Proof_QP_bcd}
As illustrated in Figs. \ref{Fig_Locations} (b)-(d), the remaining derivations correspond to lines \ref{Alg_bcd}-\ref{Alg_bcd_Result}. They examine all possible cases that the optimality occurs, corresponding to the boundary as determined by additional inequality constraint(s). It is worth noting that such a boundary can be either ``edge'' or ``vertex'', depending on its dimension. Specifically, at line \ref{Alg_For_All}, we examine all the possible, ($2^\kappa-1$)'s subsets/cases $\calI_k\subseteq\calI$. The feasibility will be checked during the process, and each one is associated with the ``augmented'' equality-constrained QP problem with respect to $(\tilde A,\tilde\bfb)$ at line \ref{Alg_bcd_Augmented_Parameters}, where $\tilde C\in\real^{\itk\times\itn}$, $\tilde A\in\real^{(\itm+\itk)\times\itn}$, and $\tilde\bfb\in\real^{\itm+\itk}$.

In the very beginning, line \ref{Alg_bcd_Condition_Consistency} filters out the inconsistent case(s): augmented QP(s) imposed by inconsistent equality constraint(s), or not an ``edge/vertex'' of the polyhedron $\calF$ as shown in Fig. \ref{Fig_Locations}. Obviously, the optimality does not occur in such case(s). The next stage, the if-else setting at lines \ref{Alg_bcd_Condition_Full_Rank_Singular}-\ref{Alg_bcd_Condition_Full_Rank_End}, filters in only the cases of full-rank $\tilde A$, according to the formulation in (\ref{Problem_Equality_QP}); while the other(s) has/have already been considered before (this case/consideration), corresponding to the reduced, full-rank, equivalent counterpart. Also, the if-else setting distinguishes between augmented systems of many solutions (commonly an underdetermined system) and of the only one. Specifically, lines \ref{Alg_bcd_Condition_Full_Rank_Singular}-\ref{Alg_bcd_Condition_Full_Rank_Singular_End} mimic the processes at lines \ref{Alg_A_V_2}-\ref{Alg_Switch_Element} but with respect to the ``augmented'' equality constraints $(\tilde A,\tilde\bfb)$, where $\tilde V_2\in\real^{\itn\times(\itn-\itm-\itk)}$ is an orthogonal matrix as computed at line \ref{Alg_bcd_A_V_2}; while lines \ref{Alg_bcd_Condition_Full_Column_Rank}-\ref{Alg_bcd_Condition_Full_Column_Rank_End} consider the full-column-rank $\tilde A$, and the unique feasible point/optimum candidate $\hat\bfx$ (computed at line \ref{Alg_bcd_Condition_Full_Column_Rank_Check_x}), which is associated with some $\calI_j$ and the only level set (corresponding to the optimal value for this augmented system). Note that:
\begin{enumerate-i-a}
\item At lines \ref{Alg_Inequality Equality_QP_Optimum_Condition_Inequality} and \ref{Alg_bcd_Condition_Full_Column_Rank_Feasibility}, it is sufficient to check all those inequality constraints that are not ``augmented/involved'' in $(\tilde A,\tilde\bfb)$ corresponding to an $\calI_j$. This is beneficial from a computational perspective, because excessive computations are avoided.
\item If the conditions at line \ref{Alg_Inequality Equality_QP_Optimum_Condition_Inequality} are all satisfied, then lines \ref{Alg_bcd_Condition_Full_Rank_Singular_Optimal_Value} and \ref{Alg_bcd_Condition_Full_Rank_Singular_Inclusion} gather the optimality candidate, as well as the involved/associated inequality constraints, in $\bar\calL\times\bar\calX\times\bar\calI$. The last element provides the information on which ``edge'' the unique optimum (or the particular solution of optima) resides, with respect to the augmented equality-constrained QP. A geometric interpretation of this \textit{optimality candidate} is illustrated in Fig. \ref{Fig_Candidates}, with more details in Remark \ref{Rem_Candidates}. On the other hand (when unsatisfied at line \ref{Alg_Inequality Equality_QP_Optimum_Condition_Inequality}), according to Theorem \ref{Thm_Equality_Constrained_QP}, the consideration instead goes to one of the following:
    \begin{enumerate}
    \item\label{App_Proof_QP_bcd_LP} An unconstrained LP, as in \ref{Thm_Equality_Constrained_QP_Condition_Unconstrained_LP}) of Theorem \ref{Thm_Equality_Constrained_QP}, the optimality can be determined in another case, say $\calI_{\bar j}\supset\calI_j$, if the optimality is actually associated with the $\calI_j$ in this case. In other words, given the unbounded LP, the optimal value only happens by imposing further (inequality) constraints, and will be at the ``terminal(s)''. Accordingly, the associated optimum/optima are further termed \textit{terminal}.
    \item\label{App_Proof_QP_bcd_Constant} A constant problem, as in \ref{Thm_Equality_Constrained_QP_Condition_Constant}) of Theorem \ref{Thm_Equality_Constrained_QP}. Obviously, its optimal value is just the only level set value. However, instead of checking whether the domain intersects $\calF$, a more efficient way is to reserve this case to another one in the for-all environment, since there are only two possibilities: one is that this case will never be further constrained, which obviously does not solve QP Problem (\ref{Problem_QP}); regarding the other, the optimal value still remains and \textit{also} occurs at the ``terminal optimum/optima'', which is/are included at line \ref{Alg_bcd_Condition_Full_Rank_Singular_Inclusion} and later compared with all the other candidate(s) at line \ref{Alg_bcd_Result}. This case is demonstrated by Example \ref{Ex_QP_Singular} (Sec. \ref{Sec_Ex}).
    \item An (equality-constrained) QP but unbounded below, as determined by {\ref{Thm_Equality_Constrained_QP_Unconstrained_QP_Finite_Condition})} of Theorem \ref{Thm_Equality_Constrained_QP}. This derivation is similar to {\ref{App_Proof_QP_bcd_LP})} above and thus omitted.
    \item\label{App_Proof_QP_bcd_QP_bounded} An (equality-constrained) QP of finite optimal value, with its unique optimum, or the particular solution of optima, outside of the polyhedron $\calF$; in other words, not satisfying the other inequality constraint(s) that is/are not involved in this augmented system. Notably, there is still only one possibility that the optimality of QP Problem (\ref{Problem_QP}) is associated with this case, which corresponds to the nonuniqueness of optima. Specifically, this possibility is owing to the further imposed equality constraint(s), or the boundary of inequality constraint(s) from $\calI\backslash\calI_j$, on the degree of freedom of the optima with respect to the augmented QP under equality constraints and additional ones from $\calI_j$. This will be elsewhere included/considered by virtue of the for-all environment in Algorithm \ref{Alg_QP}.
    \end{enumerate}
\item Regarding the full-column-rank $\tilde A$, at lines \ref{Alg_bcd_Condition_Full_Column_Rank}-\ref{Alg_bcd_Condition_Full_Column_Rank_End}, the only feasible point/optimum with respect to the augmented equality-constrained QP is computed a priori at line \ref{Alg_bcd_Condition_Full_Column_Rank_Check_x}. If this point is also feasible with respect to the original/general QP Problem (\ref{Problem_QP}), as determined at line \ref{Alg_bcd_Condition_Full_Column_Rank_Feasibility}, then the corresponding, only, feasible/optimal value is later computed at line \ref{Alg_bcd_Condition_Full_Column_Rank_Optimal_Value}, and included as an optimality candidate at line \ref{Alg_bcd_Condition_Full_Column_Rank_Inclusion} (the overall comparison in the end of this algorithm).
\end{enumerate-i-a}

Finally, at line \ref{Alg_bcd_Result}, the result concludes by comparing all the candidates in $\bar\calL$ of finite elements. As a matter of fact, the cardinality of $\bar\calL$ is at most $2^\kappa$, which is detailed later in Remark \ref{Rem_Candidates_Cardinality}. The minimum/optimal value $\tilde l^*$ is associated with the corresponding \textit{terminal optimum/optima} in $\bar\calX$, which can be either the unique optimum, or the particular/terminal solution of nonunique optima, located on the ``edge(s)'' or ``vertex/vertices'' as determined by $\tilde\calI^*$.

\begin{Remark}
In the end of proof/algorithm, such an optimality pair $(\tilde l^*,\bfx^{\tilde *},\tilde\calI^*)$ can be nonunique. A more in-depth interpretation is that the optimality can occur at the different/nonunique optima, or different ``edges/vertices''. The information regarding the latter is completely gathered by means of this algorithm. As for the former, if the optimum is unique with respect to all possible $\tilde\calI^*$('s), then this information is also complete. On the other hand, if the optima (with respect to some $\tilde\calI^*$) are nonunique, then, to remain focused of this presentation, preliminarily a subset is formulated using Algorithm \ref{Alg_QP}. This subset includes all the \textit{terminal optima}. That is, Algorithm \ref{Alg_QP} preliminarily gives closed-form solutions for all the optimum/optima at the terminals; while the intermediate points require specific attention due to the flexibility of inequality constraints. Although this indicates an interesting further research direction, Example \ref{Ex_QP_Singular} in Sec. \ref{Sec_Ex} demonstrates a potential of the proposed results for such a comprehensive generalization. That is, a closed-form solution to the considered QP, including an explicit representation for all the (nonunique) optima.
\label{Rem_Terminal_Optima}
\end{Remark}

\begin{Remark}
Fig. \ref{Fig_Candidates} demonstrates a fundamental concept ``optimality candidates'' for the design of Algorithm \ref{Alg_QP}, in terms of the case where all $\bfx^{\bar *}\not\in\calF$. It is worth noting that, according to Eq. (\ref{Optimum_Equality_Constrained_QP_Unconstrained_QP}), if there exists an $\bfvarepsilon^{\bar *}$ such that the corresponding $\bfx^{\bar *}\in\calF$, then $\tilde l^*=\bar l^*$, and follow similar discussions in the remaining of this remark. Regarding this planar example of Fig. \ref{Fig_Candidates}, the optimality occurs at (d) $\bfx^{\tilde *}$, and there are three representative cases/points (a)-(c) that are considered in the for-all design at lines \ref{Alg_For_All}-\ref{Alg_For_All_End} of Algorithm \ref{Alg_QP}:
\begin{enumerate-a}
\item The corresponding augmented pair $(\tilde A,\tilde\bfb)$ satisfies all the conditions at lines \ref{Alg_bcd_Condition_Full_Rank_Singular} and \ref{Alg_Inequality Equality_QP_Optimum_Condition_Inequality}, except that its unique optimum is outside of the feasible set $\calF$. Thus, this point is not an optimality candidate.
\item The pair $(\tilde A, \tilde\bfb)$ satisfies all the conditions at lines \ref{Alg_bcd_Condition_Full_Rank_Singular} and \ref{Alg_Inequality Equality_QP_Optimum_Condition_Inequality}, where its unique optimum is located on an edge of $\calF$. Hence, this optimality candidate $\tilde\bfx_p^{\bar *}=\tilde\bfx^{\bar *}$ is included at line \ref{Alg_bcd_Condition_Full_Rank_Singular_Inclusion}, for the overall comparison later at line \ref{Alg_bcd_Result}.
\item Since $\tilde A$ is of full column rank and its associated unique feasible point (optimum) $\hat\bfx$ is also feasible with regard to the other inequality constraints, $\hat\bfx$ is a candidate and also a vertex of $\calF$, which is included at line \ref{Alg_bcd_Condition_Full_Column_Rank_Inclusion} for the overall comparison (line \ref{Alg_bcd_Result}).
\end{enumerate-a}
However, the optimality of this demonstration occurs at the vertex/candidate (d) $\bfx^{\tilde *}$, whose arguments similarly follow the above (c). After comparing the associated, finitely many optimal values among all the candidates (line \ref{Alg_bcd_Result}), including the ones at (b)-(d), the corresponding unique optimum $\bfx^{\tilde *}$, optimal (level set) value $\tilde l^*$, and information for the location of $\bfx^{\tilde *}$ ($\tilde\calI^*$) are exactly determined. Note that the complete and explicit formulation of the QP example in Fig. \ref{Fig_Candidates} is not required for the discussions in this remark, particularly regarding the inequality constraints, since it is compactly two-dimensional for ease of reading.
\label{Rem_Candidates}
\end{Remark}

\begin{Remark}
This is to determine the largest possible cardinality of $\bar\calL$. To start with, note that there are $2^\kappa$'s different subsets of $\calI$. The case of empty set is considered in lines \ref{Alg_a_optimum_particular}-\ref{Alg_a_QP_Condition_End} to check the feasibility of $\bfx_p^{\bar *}$, which is associated with $\bar\calL=\emptyset$ (as in the designed initial condition); while its possible freedom $(V_2\bfvarepsilon^{\bar *})$ is designed to be checked within the for-all environment (lines \ref{Alg_For_All}-\ref{Alg_For_All_End}), if the solutions to QP Problems (\ref{Problem_Equality_QP}) and (\ref{Problem_QP}) intersect, such as the demonstrations (a) and (b) in Fig. \ref{Fig_Locations}. On the other hand, the for-all design/coverage includes all the other subsets (note that $k\ge 1$ at line \ref{Alg_For_All}), which is associated with, at most, $(2^\kappa-1)$'s elements for $\bar\calL$.
\label{Rem_Candidates_Cardinality}
\end{Remark}

\begin{Remark}
An extreme example considers the singleton $\calF\ne\emptyset$ (noting that the subset $\emptyset$ contains no point), and this remark explains how Algorithm \ref{Alg_QP} solves it. Similarly, divide the derivations into whether $\{\bfx_p^{\bar *}\}=\calF$. Regarding the equivalent/simpler case, then the solving process can be more efficiently completed without excessively going into the for-all environment (lines \ref{Alg_For_All}-\ref{Alg_For_All_End}), while returning and ensuring the correct optimality. On the other hand, the process at lines \ref{Alg_bcd_Condition_Full_Rank_Singular}-\ref{Alg_bcd_Condition_Full_Rank_Singular_End} yields the optimality, only if $\tilde\bfx_p^{\bar *}\in\calF$ for some $\calI_j$; while that at \ref{Alg_bcd_Condition_Full_Column_Rank}-\ref{Alg_bcd_Condition_Full_Column_Rank_End} ``exhaustively'' includes this case, a vertex of $\calF$. Although in the former (two) cases, this algorithm only concentrates on the particular solution $\bfx_p^{\bar *}/\tilde\bfx_p^{\bar *}$; nonetheless, this example also demonstrates that the consideration on the freedom $(V_2\bfvarepsilon^{\bar *}/\tilde V_2\tilde\bfvarepsilon^{\bar *})$ can be more than necessary, where, for brevity, the notation $\tilde\bfvarepsilon^{\bar *}$ similarly follows $\bfvarepsilon^{\bar *}$ but with respect to the augmented pair $(\tilde A,\tilde\bfb)$. \textit{As a comparison}, it is not straightforward to apply most of the existing QP solvers \cite{BoVa:04,NoWr:06} to this singleton example, which are commonly/classically designed from the differential perspective while initiated by (knowing) a feasible point, such as the popular ASM \cite{Lu(Ye):03(16)}.
\label{Rem_QP_Singleton_F}
\end{Remark}

\section{Proof of Theorem \ref{Thm_Extended_QP} (Extended QP)}
\label{App_Proof_Extended_QP}

\noindent\ref{Thm_Extended_QP_Preimage}) By \ref{Thm_Unconstrained_QP_all_preimage}) of Theorem \ref{Thm_Unconstrained_QP}, rewrite equivalently the CCQF at any level set value of $\check l$ as the following constrained convex quadratic equation (CCQE):
\beq
&&(1/2)\bfx^TP\bfx+\bfq^T\bfx+s-\check l=0\nonumber\\
&&\mbox{subject to both}~\bfx,~\bfq\in \calR(P).
\label{Thm_Extended_QP_CCQE}
\eeq
According to Lemma \ref{Lem_Solutions}, a necessary condition that renders the above CCQE in (\ref{Thm_Extended_QP_CCQE}) solvable is
\begin{subnumcases}{\label{Solvability_Cond_Thm_Extended_QP}}
\bfq^TP^{-1}\bfq-2s+2\check l\ge 0,      & if rank$(P)=\itn$;\label{Solvability_Cond_Thm_Extended_QP_n}\\
\bfq^TP^\dagger\bfq-2s+2\check l\ge 0~\mbox{and}~\bfx\in \calR(P), & if rank$(P)<\itn$.\label{Solvability_Cond_Thm_Extended_QP_r}
\end{subnumcases}
Note that the constraint ``$\bfq\in\calR(P)$'' has been (implicitly) taken into consideration. The remaining part of this proof starts with the derivations for the sufficiency of Condition (\ref{Solvability_Cond_Thm_Extended_QP}). Accordingly, divide the derivations into whether the Hessian matrix is of full rank or not, as \ref{Thm_Extended_QP_Preimage}-1) and \ref{Thm_Extended_QP_Preimage}-2) below.\vspace{0.16cm}

\noindent\ref{Thm_Extended_QP_Preimage}-1) rank$(P)=\itn$. In this case, we have i) $P$ is nonsingular, ii) $\calR(P)=\real^\itn$, and thus iii) the constraint ``both $\bfx,~\bfq\in\calR(P)$'' is automatically lifted, which ensures the sufficiency of Condition (\ref{Solvability_Cond_Thm_Extended_QP_n}). By virtue of Eq./Parameterization (\ref{Solution_n}) under Condition (\ref{Solvability_Cond_Thm_Extended_QP_n}), the result thus follows.\vspace{0.16cm}

\noindent\ref{Thm_Extended_QP_Preimage}-2) rank$(P)=\itr<\itn$. Let the SVD of $P$ be given by
\beq
P&=&[\tilde U_1~\tilde U_2]\left[ \begin{array}{cc} \tilde\Sigma_1~~ & 0\\ 0~~ & 0\end{array}\right]\left[ \begin{array}{c} \tilde U_1^T\\\tilde U_2^T\end{array}\right]\nonumber\\
&=&\tilde U_1\tilde\Sigma_1\tilde U_1^T,\label{SVD_P}
\eeq
where $\tilde U_1\in\real^{{\it n}\times {\it r}}$, $\tilde U_2\in\real^{{\it n}\times ({\it n-r})}$, and $\tilde\Sigma_1\in\real^{{\it r}\times {\it r}}$. For conciseness, the properties of this SVD (\ref{SVD_P}) for derivations afterward are similar to, and thus referenced at, \cite[Eqs.~{(35)-(37)}]{LiLiHs:20}. Likewise, define $P^\dagger\coloneqq\tilde U_1\tilde\Sigma_1^{-1}\tilde U_1^T$ and $P^{\dagger/2}\coloneqq\tilde U_1\tilde\Sigma_1^{-1/2}\tilde U_1^T$. As a result, the CCQE in (\ref{Thm_Extended_QP_CCQE}) is equivalent to the following CQE:
\beq
&&(1/2)\bfx^T\tilde U_1\tilde\Sigma_1\tilde U_1^T\bfx+\bfq^T(\tilde U_1\tilde U_1^T+\tilde U_2\tilde U_2^T)\bfx+s-\check l=0\nonumber\\
&&\Leftrightarrow (1/2)\check\bfx^T\tilde\Sigma_1\check\bfx+\hat\bfq^T\check\bfx+s-\breve F_{\bfx}=0,
\label{Thm_Extended_QP_equi_CQE}
\eeq
where $\tilde U_1\tilde U_1^T+\tilde U_2\tilde U_2^T=I_\itn$ is inserted purposely, and $\check\bfx=\tilde U_1^T\bfx\in\real^\itr$ (resp., $\hat\bfq=\tilde U_1^T\bfq\in\real^\itr$) is the coordinate vector of $\bfx\in \calR(P)$ (resp., $\bfq\in \calR(P)$) with respect to the $\tilde U_1$-basis. Note that the original CCQF has been equivalently transformed into the CQE (\ref{Thm_Extended_QP_equi_CQE}) with the full-rank Hessian matrix $\tilde\Sigma_1$, and the dimension is reduced from $\itn$ to $\itr$ (because of the constraint ``both $\bfx,~\bfq\in\calR(P)$''). In addition, the solvability of the CQE (\ref{Thm_Extended_QP_equi_CQE}) has already been ensured. To see this, by \ref{Lem_Solutions_n}) of Lemma \ref{Lem_Solutions}, the CQE (\ref{Thm_Extended_QP_equi_CQE}) is solvable, iff
\beq
&&(1/4)\hat\bfq^T(\tilde\Sigma_1/2)^{-1}\hat\bfq\ge s-\check l\nonumber\\
&&\Leftrightarrow \bfq^T\tilde U_1\tilde\Sigma_1^{-1}\tilde U_1^T\bfq-2s+2\check l\ge 0.\nonumber
\eeq
This has already been given/guaranteed according to Condition/Eq. (\ref{Solvability_Cond_Thm_Extended_QP_r}), and also ensures its sufficiency that renders the above CCQE (\ref{Thm_Extended_QP_CCQE}) solvable.

Therefore, by virtue of Parameterization/Eq. (\ref{Solution_n}), the solution set to the CQE in (\ref{Thm_Extended_QP_equi_CQE}), where $\hat\Sigma_1$ is of full rank ($\itr$), is given as
\beq
\check\bfx=-\tilde\Sigma_1^{-1}\hat\bfq+\sqrt{{\hat\bfq^T\tilde\Sigma_1^{-1}\hat\bfq/2}-s+\check l}\cdot ({\tilde\Sigma_1/2})^{-{1/2}}\hat\bfv,~
\eeq
where $\hat\bfv\in\real^\itr$ and $\Vert\hat\bfv\Vert=1$. In terms of the original coordinate, finally we obtain all the constrained solutions (namely, elements):
\beq
\bfx&{}={}&\tilde U_1\check\bfx\nonumber\\
&{}={}&-\tilde U_1\tilde\Sigma_1^{-1}\tilde U_1^T\bfq+\sqrt{\bfq^T\tilde U_1\tilde\Sigma_1^{-1}\tilde U_1^T\bfq/2-s+\check l}\cdot\sqrt{2}\cdot\tilde U_1\tilde\Sigma_1^{-1/2}\cdot\tilde U_1^T\tilde U_1\cdot \hat\bfv\nonumber\\
&{}={}&\mbox{Eq. (\ref{Thm_Extended_QP_Solutions})},\nonumber
\eeq
where $\tilde U_1^T\tilde U_1=I_\itr$ is inserted purposely, $\check\bfrho=\tilde U_1\hat\bfv\in \calR(P)$, and $\Vert\check\bfrho\Vert=\Vert\tilde U_1\hat\bfv\Vert=\Vert\hat\bfv\Vert=1$. Remarkably, the constraint $\bfx\in \calR(P)$ is satisfied since the two vectors in Eq. (\ref{Thm_Extended_QP_Solutions}), $P^\dagger\bfq$ and $P^{\dagger/2}\check\bfrho$, reside in $\calR(P)$.\vspace{0.16cm}

\noindent\ref{Thm_Extended_QP_Optimality}) The optimal value is finite, whose arguments similarly follow the proof for \ref{Thm_Unconstrained_QP_finite_condition}) of Theorem \ref{Thm_Unconstrained_QP}, and is mainly due to the constraint ``$\bfq\in\calR(P)$''. Moreover, according to the above derivations for \ref{Thm_Extended_QP_Preimage}), the necessary condition in Eq. (\ref{Solvability_Cond_Thm_Extended_QP}) is also sufficient to render the CCQE (\ref{Thm_Extended_QP_CCQE}) solvable. Thus, the optimal value $\check l^*$ readily follows from the Condition (\ref{Solvability_Cond_Thm_Extended_QP}). Substituting $\check l^*$ into Eq. (\ref{Thm_Extended_QP_Solutions}) yields the associated unique optimum $\bfx^{\hat *}$, which completes the proof.

\begin{Remark}
A geometric interpretation of the preimage, as addressed in \ref{Thm_Extended_QP_Preimage}) of Theorem \ref{Thm_Extended_QP}, can be similarly inferred from Fig. \ref{Fig_Layers}, where the solution freedom residing in $\calN(P)$ is not available/effective. This is essentially due to the constraint ``$\bfx\in\calR(P)$'' in Problem (\ref{Problem_Extended_QP}), while the other ``$\bfq\in\calR(P)$'' mainly renders this optimization problem with a finite optimal value. From a different viewpoint, this remark is also consistent with Eq. (\ref{Solution_k_R(M)}), whose derivations are similar and hence omitted for conciseness.
\end{Remark}

\section{Tables for all $\calI_j$'s in Examples \ref{Ex_QP_Nonsingular} and \ref{Ex_QP_Singular} (Sec. \ref{Sec_Ex})}
\label{App_Tables}

This appendix is designed to facilitate comparisons between Examples \ref{Ex_QP_Nonsingular} and \ref{Ex_QP_Singular}. Table \ref{Table_Ex_QP_Nonsingular} (resp., \ref{Table_Ex_QP_Singular}) is associated with Example \ref{Ex_QP_Nonsingular} (resp., \ref{Ex_QP_Singular}), which demonstrates the solving process on a \textit{positive definite Hessian matrix} given in \cite{Lu(Ye):03(16)} (resp., \textit{positive semidefinite and singular Hessian matrix}), to summarize the details/steps when examining each case $\calI_j$ according to Algorithm \ref{Alg_QP}, specifically, the for-all environment at lines \ref{Alg_For_All}-\ref{Alg_For_All_End}. Note that the case of subset/empty set ``$\emptyset$'' (in other words, no inequality constraint imposed) has been more efficiently considered a priori, in the beginning of the two examples (lines \ref{Alg_a_QP_Condition}-\ref{Alg_a_QP_Condition_End} in Algorithm \ref{Alg_QP}). Moreover, $\calI_7=\{1,2,3\}$ imposed by all the three inequality constraints is excluded in both tables/examples, for compactness, because it is obviously not associated with the optimality ($\tilde\bfb\not\in\calR(\tilde A)$).

\begin{table}[htbp]
\centering
\begin{threeparttable}[tbhp]
\caption{A summary of significant parameters and values for all $\calI_j$'s in Example \ref{Ex_QP_Nonsingular}}
\renewcommand{\arraystretch}{1.7}
\begin{tabular}{l l l l l l l l l}
\hline
$\calI_j~(\tilde\bfb\in\calR(\tilde A))$ & $\tilde A$ & $\tilde\bfb$ & $\tilde V_2$ & $\tilde V_2^TP\tilde V_2$ & Categorization & $\tilde \bfx^{\bar *}/\hat\bfx$? & Candidate? & $\tilde{\bar l}^*/\hat l$\hspace{0.05cm}? \\ \hline \vspace{-0.8cm} \\
$\calI_1=\{1\}$ & \hspace{-0.07cm}$\left[ \renewcommand{\arraystretch}{1} \begin{array}{cc} 1~~ & 1 \end{array}\right]$ & 4 & \hspace{-0.1cm}$\left[ \renewcommand{\arraystretch}{1} \begin{array}{cc} {1\over\sqrt{2}} & ~{-1\over\sqrt{2}} \end{array}\right]^T$ & 2 & QP & $\tilde\bfx^{\bar *}=\hspace{-0.1cm}\left[ \renewcommand{\arraystretch}{1} \begin{array}{cc} 1.5 & ~~2.5 \end{array}\right]^T$ & yes & $\tilde{\bar l}^*=-28.5$ \vspace{0.1cm} \\
$\calI_2=\{2\}$ & \hspace{-0.07cm}$\left[ \renewcommand{\arraystretch}{1} \begin{array}{cc} -1~~ & 0 \end{array}\right]$ & 0 & \hspace{-0.1cm}$\left[ \renewcommand{\arraystretch}{1} \begin{array}{cc} 0~~ & 1 \end{array}\right]^T$ & 2 & QP & $\tilde\bfx^{\bar *}=\hspace{-0.1cm}\left[ \renewcommand{\arraystretch}{1} \begin{array}{cc} 0~~ & 5 \end{array}\right]^T$ & no & n/a \vspace{0.1cm} \\
$\calI_3=\{3\}$ & \hspace{-0.07cm}$\left[ \renewcommand{\arraystretch}{1} \begin{array}{cc} 0~~ & -1 \end{array}\right]$ & 0 & \hspace{-0.1cm}$\left[ \renewcommand{\arraystretch}{1} \begin{array}{cc} -1~~ & 0 \end{array}\right]^T$ & 4 & QP & $\tilde\bfx^{\bar *}=\hspace{-0.1cm}\left[ \renewcommand{\arraystretch}{1} \begin{array}{cc} 3~~ & 0 \end{array}\right]^T$ & yes & $\tilde{\bar l}^*=-18$ \vspace{0.1cm} \\
$\calI_4=\{1,2\}$ & \hspace{-0.1cm}$\left[ \renewcommand{\arraystretch}{1} \begin{array}{cc} 1 & ~~1 \\ -1 & ~~0 \end{array}\right]$ & \hspace{-0.1cm}$\left[ \renewcommand{\arraystretch}{1} \begin{array}{c} 4 \\ 0 \end{array}\right]$ & \hspace{-0.05cm}n/a & n/a & vertex & $\hat\bfx=\hspace{-0.1cm}\left[ \renewcommand{\arraystretch}{1} \begin{array}{c} 0 \\ 4 \end{array}\right]$ & yes & $\hat l=-24$ \vspace{0.1cm} \vspace{0.1cm} \\
$\calI_5=\{2,3\}$ & \hspace{-0.1cm}$\left[ \renewcommand{\arraystretch}{1} \begin{array}{cc} -1 & ~~0 \\ 0 & ~~{-1} \end{array}\right]$ & \hspace{-0.1cm}$\left[ \renewcommand{\arraystretch}{1} \begin{array}{c} 0 \\ 0 \end{array}\right]$ & \hspace{-0.05cm}n/a & n/a & vertex & $\hat\bfx=\hspace{-0.1cm}\left[ \renewcommand{\arraystretch}{1} \begin{array}{c} 0 \\ 0 \end{array}\right]$ & yes & $\hat l=0$ \vspace{0.1cm} \vspace{0.1cm} \\
$\calI_6=\{1,3\}$ & \hspace{-0.1cm}$\left[ \renewcommand{\arraystretch}{1} \begin{array}{cc} 1 & ~~1 \\ 0 & ~~{-1} \end{array}\right]$ & \hspace{-0.1cm}$\left[ \renewcommand{\arraystretch}{1} \begin{array}{c} 4 \\ 0 \end{array}\right]$ & \hspace{-0.05cm}n/a & n/a & vertex & $\hat\bfx=\hspace{-0.1cm}\left[ \renewcommand{\arraystretch}{1} \begin{array}{c} 4 \\ 0 \end{array}\right]$ & yes & $\hat l=-16$ \vspace{0.2cm} \\
\hline
\end{tabular}
\label{Table_Ex_QP_Nonsingular}
\end{threeparttable}
\end{table}

\begin{table}[htbp]
\centering
\begin{threeparttable}[htbp]
\caption{A summary of significant parameters and values for all $\calI_j$'s in Example \ref{Ex_QP_Singular}}
\renewcommand{\arraystretch}{1.7}
\begin{tabular}{l l l l l l l l l l}
\hline
$\calI_j~(\tilde\bfb\in\calR(\tilde A))$ & \hspace{-0.03cm}$\tilde A$ & \hspace{-0.25cm}$\tilde\bfb$ & \hspace{-0.2cm}$\tilde V_2$ & \hspace{-0.2cm}$\tilde V_2^TP\tilde V_2$ & \hspace{-0.1cm}Categorization & \hspace{-0.1cm}$\calN(P)\cap\calN(\tilde A)$ & \hspace{-0.1cm}$\tilde \bfx^{\bar *}/~\hat\bfx$? & \hspace{-0.1cm}Candidate? & \hspace{-0.1cm}$\tilde{\bar l}^*/~\hat l$\hspace{0.05cm}? \\ \hline \vspace{-0.6cm} \\
$\calI_1=\{1\}$ & \hspace{-0.1cm}$\left[ \renewcommand{\arraystretch}{1} \begin{array}{ccc} 0 & ~~0~~ & 1 \\ 1 & ~~1~~ & 0 \end{array}\right]$ & \hspace{-0.3cm}$\left[ \renewcommand{\arraystretch}{1} \begin{array}{c} 0\\ 4\end{array}\right]$~ & \hspace{-0.2cm}$\left[ \renewcommand{\arraystretch}{1} \begin{array}{c} {-1\over\sqrt{2}} \\ {1\over\sqrt{2}} \\ 0 \end{array}\right]$~ & \hspace{-0.2cm}${1\over 2}$ & \hspace{-0.1cm}QP & \hspace{-0.1cm}\{\bfzero\} & \hspace{-0.1cm}$\tilde\bfx^{\bar *}=\left[ \renewcommand{\arraystretch}{1} \begin{array}{c} 0 \\ 4 \\ b\end{array}\right]$~ & \hspace{-0.1cm}yes & \hspace{-0.1cm}$\tilde{\bar l}^*=0$ \vspace{0.1cm} \vspace{0.1cm} \\
$\calI_2=\{2\}$ & \hspace{-0.1cm}$\left[ \renewcommand{\arraystretch}{1} \begin{array}{ccc} 0 & ~~0~~ & 1 \\ -1 & ~~0~~ & 0 \end{array}\right]$ & \hspace{-0.3cm}$\left[ \renewcommand{\arraystretch}{1} \begin{array}{c} b\\ 0\end{array}\right]$ & \hspace{-0.2cm}$\left[ \renewcommand{\arraystretch}{1} \begin{array}{c} 0 \\ {-1} \\ 0 \end{array}\right]$ & \hspace{-0.2cm}0 & \hspace{-0.1cm}constant\tnote{a} & \hspace{-0.1cm}n/a & \hspace{-0.1cm}n/a & \hspace{-0.1cm}n/a\tnote{b} & \hspace{-0.1cm}n/a \vspace{0.1cm} \vspace{0.1cm} \\
$\calI_3=\{3\}$ & \hspace{-0.1cm}$\left[ \renewcommand{\arraystretch}{1} \begin{array}{ccc} 0 & ~~0~~ & 1 \\ 0 & ~~{-1}~~ & 0 \end{array}\right]$ & \hspace{-0.3cm}$\left[ \renewcommand{\arraystretch}{1} \begin{array}{c} b\\ 0\end{array}\right]$ & \hspace{-0.2cm}$\left[ \renewcommand{\arraystretch}{1} \begin{array}{c} 1 \\ 0 \\ 0 \end{array}\right]$ & \hspace{-0.2cm}1 & \hspace{-0.1cm}QP & \hspace{-0.1cm}\{\bfzero\} & \hspace{-0.1cm}$\tilde\bfx^{\bar *}=\left[ \renewcommand{\arraystretch}{1} \begin{array}{c} 0 \\ 0 \\ b\end{array}\right]$ & \hspace{-0.1cm}yes & \hspace{-0.1cm}$\tilde{\bar l}^*=0$ \vspace{0.1cm} \vspace{0.1cm} \\
$\calI_4=\{1,2\}$ & \hspace{-0.1cm}$\left[ \renewcommand{\arraystretch}{1} \begin{array}{ccc} 0 & ~~0~~ & 1 \\ 1 & ~~1~~ & 0 \\ -1 & ~~0~~ & 0 \end{array}\right]$ & \hspace{-0.3cm}$\left[ \renewcommand{\arraystretch}{1} \begin{array}{c} b\\ 4\\ 0\end{array}\right]$ & \hspace{-0.15cm}n/a & \hspace{-0.2cm}n/a & \hspace{-0.1cm}vertex & \hspace{-0.1cm}n/a & \hspace{-0.1cm}$\hat\bfx=\left[ \renewcommand{\arraystretch}{1} \begin{array}{c} 0\\ 4\\ b\end{array}\right]$ & \hspace{-0.1cm}yes & \hspace{-0.1cm}$\hat l=0$ \vspace{0.1cm} \vspace{0.1cm} \\
$\calI_5=\{2,3\}$ & \hspace{-0.1cm}$\left[ \renewcommand{\arraystretch}{1} \begin{array}{ccc} 0 & ~~0~~ & 1 \\ -1 & ~~0~~ & 0 \\ 0 & ~~{-1}~~ & 0 \end{array}\right]$ & \hspace{-0.3cm}$\left[ \renewcommand{\arraystretch}{1} \begin{array}{c} b\\ 0\\ 0\end{array}\right]$ & \hspace{-0.15cm}n/a & \hspace{-0.2cm}n/a & \hspace{-0.1cm}vertex & \hspace{-0.1cm}n/a & \hspace{-0.1cm}$\hat\bfx=\left[ \renewcommand{\arraystretch}{1} \begin{array}{c} 0\\ 0\\ b\end{array}\right]$ & \hspace{-0.1cm}yes & \hspace{-0.1cm}$\hat l=0$ \vspace{0.1cm} \vspace{0.1cm} \\
$\calI_6=\{1,3\}$ & \hspace{-0.1cm}$\left[ \renewcommand{\arraystretch}{1} \begin{array}{ccc} 0 & ~~0~~ & 1 \\ 1 & ~~1~~ & 0 \\ 0 & ~~{-1}~~ & 0 \end{array}\right]$ & \hspace{-0.3cm}$\left[ \renewcommand{\arraystretch}{1} \begin{array}{c} b\\ 4\\ 0\end{array}\right]$ & \hspace{-0.15cm}n/a & \hspace{-0.2cm}n/a & \hspace{-0.1cm}vertex & \hspace{-0.1cm}n/a & \hspace{-0.1cm}$\hat\bfx=\left[ \renewcommand{\arraystretch}{1} \begin{array}{c} 4\\ 0\\ b\end{array}\right]$ & \hspace{-0.1cm}yes & \hspace{-0.1cm}$\hat l=16$ \vspace{0.2cm} \\
\hline
\end{tabular}
\begin{tablenotes}
\item [a] $\tilde V_2^T\bfq+\tilde V_2^TP\tilde A^\dagger\tilde \bfb=0$ (not an LP), and the constant value is $s+\bfq^T\tilde A^\dagger\tilde\bfb+\tilde\bfb^T(\tilde A^\dagger)^TP\tilde A^\dagger\tilde\bfb/2=0$.
\item [b] The associated \textit{terminal optima} are explicitly formulated after further imposing constraints $(\bfc_1,d_1)$ and $(\bfc_3,d_3)$, as included in  Cases $\{\calI_4\}$ and $\{\calI_5\}$, respectively.
\end{tablenotes}
\label{Table_Ex_QP_Singular}
\end{threeparttable}
\end{table}

\vspace{4cm}

\section*{Acknowledgment}
The authors acknowledge Julian Mellor and Scott Read for their assistance in the native use of this language. This work was supported in part by the Max Planck Society, Germany; Ministry of Science and Technology, Taiwan; and Delta Electronics, Inc.

%\bibliographystyle{amsplain}
%\bibliography{2022_SDRE}

\providecommand{\bysame}{\leavevmode\hbox to3em{\hrulefill}\thinspace}
\providecommand{\MR}{\relax\ifhmode\unskip\space\fi MR }
\providecommand{\MRhref}[2]{%
  \href{http://www.ams.org/mathscinet-getitem?mr=#1}{#2}
}
\providecommand{\href}[2]{#2}

\end{document}